\documentclass[12pt]{article}
\usepackage{lmodern}
\usepackage[margin=2cm,a4paper]{geometry}
\usepackage{graphicx} %
\usepackage{amssymb}
\usepackage{amsmath}
\usepackage{xcolor}
\usepackage{multirow}

\usepackage{hyperref}
\usepackage{authblk}
\usepackage{zref-clever}
\zcsetup{cap=true}
\usepackage{enumitem}
\usepackage[thmtools-compat]{keytheorems}

\usepackage{mathrsfs}
\usepackage{tikz}
\usetikzlibrary{shapes.geometric, positioning, fit, backgrounds, calc}

\hypersetup{
    colorlinks=true,
    citecolor=blue,
    linkcolor=blue,
    filecolor=blue,      
    urlcolor=blue,
    pdftitle={Blind cop-width and balanced minors of graphs},
    pdfauthor={Hector Buffière, Rutger Campbell, Kevin Hendrey, and Sang-il Oum}
}

\newcommand{\mynote}[3]{
		\textcolor{#2}{\fbox{\bfseries\sffamily\scriptsize#1}}
		{\small\textsf{\emph{\textcolor{#2}{#3}}}}}

\newcommand{\kh}[1]{\mynote{Kevin}{blue}{#1}}
\newcommand\kevin\kh
\newcommand{\rc}[1]{\mynote{Rutger}{violet}{#1}}
\newcommand\rutger\rc

\newcommand\abs[1]{\left\lvert#1\right\rvert}

\declaretheorem[name=Theorem, numberwithin=section, refname={Theorem,Theorems}, Refname={Theorem,Theorems}]{thm}

\declaretheorem[name=Lemma, sibling=thm, refname={Lemma,Lemmas}, Refname={Lemma,Lemmas}]{lem}
\declaretheorem[name=Proposition, sibling=thm, refname={Proposition,Propositions}, Refname={Proposition,Propositions}]{prop}
\declaretheorem[name=Corollary, sibling=thm, refname={Corollary,Corollaries}, Refname={Corollary,Corollaries}]{cor}
\declaretheorem[name=Conjecture, sibling=thm, refname={Conjecture,Conjectures}, Refname={Conjecture,Conjectures}]{conj}
\declaretheorem[name=Definition, sibling=thm, refname={Definition,Definitions}, Refname={Definition,Definitions}]{definition}

\declaretheorem[name=Problem, sibling=thm, refname={Problem,Problems}, Refname={Problem,Problems}]{problem}

\newcommand\pw{\operatorname{pw}}
\newcommand\tw{\operatorname{tw}}
\newcommand\copwidth{\operatorname{copwidth}}
\newcommand\degeneracy{\operatorname{degeneracy}}
\newcommand\bfw{\operatorname{bfw}}
\newcommand\had{\operatorname{h}}
\newcommand\tbcw{\operatorname{tbcw}}
\newcommand\fw{\operatorname{fw}}

\newcommand\IN{\operatorname{in}}
\newcommand\HH{\operatorname{hunt}}
\newcommand\fc{\operatorname{\textsf{FC}}}
\newcommand\pc{\operatorname{\textsf{PC}}}
\AtEndEnvironment{ex}{\hfill$\triangle$}
\let\geq=\geqslant
\let\leq=\leqslant
\let\ge=\geqslant 
\let\le=\leqslant

\newcommand\bcw{\operatorname{bcw}}
\newcommand\minor\preccurlyeq

\title{Blind cop-width and balanced minors of graphs}
\author[1,2]{Hector Buffière\thanks{Supported by the Institute for Basic Science (IBS-R029-C1).}%
\thanks{Supported by 
the European Research Council (ERC) under the European
Union's Horizon 2020 research and innovation programme
(grant agreement No 810115 - Dynasnet).}}
\author[3]{Rutger Campbell\textsuperscript{\textasteriskcentered}%
\thanks{Supported by the National Research Foundation of Korea (NRF) grant funded by the Ministry of Science and ICT (No.~RS-2025-00563533).}}
\author[4]{Kevin Hendrey\textsuperscript{\textasteriskcentered}\thanks{Supported by
the Australian Research Council.}}
\author[5,6]{Sang-il~Oum\textsuperscript{\textasteriskcentered}}
\affil[1]{IRIF, Université Paris Cité, Paris, France}
\affil[2]{CAMS, École des Hautes Études en Sciences Sociales, Paris, France}
\affil[3]{School of Computing, KAIST, Daejeon, South~Korea}
\affil[4]{School of Mathematics, Monash University, Melbourne, Australia}
\affil[5]{Discrete Mathematics Group, Institute for Basic Science (IBS), Daejeon,~South~Korea}
\affil[6]{Department of Mathematical Sciences, KAIST, Daejeon, South~Korea}
\affil[ ]{\small \textit{Email addresses:} \texttt{buffiere@irif.fr},
\texttt{rutger@kaist.ac.kr},
\texttt{Kevin.Hendrey1@monash.edu},
\texttt{sangil@ibs.re.kr}}

\begin{document}
\maketitle
\begin{abstract}
	We investigate a pursuit-evasion game on an undirected graph in which a robber, moving at a fixed constant speed, attempts to evade a team of cops who are blind to the robber's location and can quickly travel between any pair of vertices in the graph. The blind cop-width is the minimum number of cops needed to catch the robber on a given graph. We link it with other known graph parameters
	defined in terms of pursuit-evasion games, and show a new lower bound with
	respect to treewidth. The proof introduces the notion of balanced
	minors, where all branch sets of a minor model have equal size.
\end{abstract}

\section{Introduction}

Pursuit-evasion games on graphs have been studied both from a game-theoretic
point of view and through their connections with structural graph theory. Cops and robber games have been identified as alternative definitions of
pathwidth and treewidth in the early days of these structural parameters. This has led them to play a crucial role in the recent progress of structural graph
theory through the graph minors project of Robertson and Seymour, and
has been used in the development of many parameterized algorithms.

In the version of the game that is used to define treewidth, a robber occupies a vertex of a graph and
tries to escape a fixed number of cops. 
At each turn, some cops use helicopters
to fly to new vertices.
Right before the helicopters land, the robber can see where the helicopters are about to land, 
and move to a new vertex by a path avoiding any vertex currently occupied by cops that did not fly. 
If the robber cannot reach an unoccupied position in the next round, the cops catch the robber and
win. In this setting, the minimum number of cops needed to be certain to catch
the robber is equal to the treewidth of the graph plus one~\cite{ST93}. If one
makes the cops unable to see the position of the robber, it is equal
to the pathwidth plus one~\cite{KP85}. 
Instead of allowing a robber to move along a path of any length, 
we may also restrict the maximum speed~$r$ of the robber.
For each integer $r$, the minimum  number of cops needed to catch the robber is then
called the radius-$r$ copwidth in \cite{T23}, where it is
shown that it is also
related to several structural graph properties.

We consider a variant called the \emph{blind cop-width game of radius $r$}, following the presentation of Toru\'nczyk~\cite{T23}.
This game imposes both conditions simultaneously, meaning that the robber is both slow and invisible. 
Initially, a robber can be hidden at any vertex. 
At each turn, some cops use helicopters to fly to new vertices.
The cops win this game if they end up catching the robber no matter where he was hidden. 
During this whole process, the cops cannot see the location of the robber.
The \emph{radius-$r$ blind cop-width} of a graph $G$, denoted by $\bcw_r(G)$, is the minimum number of cops to catch the robber in this game.

Unlike in the visible robber setting,
all radius-$r$ blind cop-width
parameters for finite~$r$ are functionally equivalent, see \zcref{sec:radius}.
In particular, for every class of graphs $\mathscr C$ and two integers~$r$ and $s$, $\bcw_r(\mathscr C)$ 
is finite if and only if $\bcw_s(\mathscr C)$ is finite
(where $\bcw_r(\mathscr{C}) = \sup_{G\in \mathscr{C}}\bcw_r(G)$). Therefore,
we focus our study on the radius-$1$ case, and derive similar results for
every other radius.
It turns out that the radius-$1$ blind cop-width is equal to the \emph{inspection number}, introduced earlier by Bernshteyn and Lee~\cite{BL22}, which we will show in \zcref{prop:bcw_in}.

Blind cop-width also shares some relations with other graph parameters defined via similar
pursuit-evasion games, which we discuss in \zcref{connections}. The first one is the zero-visibility cop number $c_0()$ \cite{DDTY15, DDTY15b, XYZZ19, XYZ22, XYZ24}, and the second one is the hunting number $\HH()$ \cite{AFGP16, BG19, DFGN23, DFGN25}. These relations allow us to get some new results about these parameters as corollaries of our results on blind cop-width.

Toru{\'n}czyk~\cite{T23} generalized cops and robber games to define \emph{flip-width}. 
Flip-width is a new graph parameter attracting a lot of attention because of its potential to be a good dense analog of the concept called `bounded expansion' of Ne{\v{s}}et{\v{r}}il and Ossona de Mendez~\cite{NO12}.  
He similarly introduced a blind variant, called \emph{blind flip-width}.
Motivated by the fact that the infinite-radius blind flip-width is functionally equivalent to linear clique-width,
Toru\'nczyk conjectured the following.
\begin{conj}[note={Toru\'nczyk~{\cite[Conjecture XI.12]{T23}}},store=torunczyk]\label{conj:torunczyk}
	A class $\mathscr{C}$ of graphs has bounded blind flip-width 
	if and only if it has bounded linear clique-width.
\end{conj}
We disprove \zcref{conj:torunczyk}. 
Indeed, we show that there is a counterexample consisting of subcubic trees.
We remark that the construction we use to refute the conjecture
had already been rediscovered several times in various equivalent forms in
independent literature about cop numbers, including Tayu and Ueno~\cite{TU15},
Bernshteyn and Lee~\cite{BL22}, 
Dereniowski, Dyer, Tifenbach, and Yang~\cite{DDTY15}, 
and Dissaux, Fioravantes, Galhawat, and Nisse~\cite{DFGN23}.

\getkeytheorem{disprove}
The proof is done by using the fact that every complete binary tree 
admits a subdivision where three cops have a winning strategy.
By generalizing this argument, we prove the following.
\getkeytheorem{treewidthsub}

\getkeytheorem{huntsub}

The \emph{topological blind cop-width} of a graph~$G$, denoted by $\tbcw(G)$, is the minimum $k$ such that every subdivision of~$G$ has a subdivision with radius-$1$ blind cop-width at most $k$.
Bernshteyn and Lee~\cite[Proposition 23]{BL22} showed that $\tbcw(G)\le \abs{V(G)}+2$, which is now a corollary of \zcref{thm:tree-sub}.
Bernshteyn and Lee~\cite{BL22} also proposed four open problems in their paper.
Three of their open problems (Problems 66, 68, and 69) can be rephrased as follows.

\begin{restatable}[Bernshteyn and Lee~{\cite[Open Problem 66]{BL22}}]{problem}{pkn}\label{p66}
   Determine $\tbcw(K_n)$ and $\tbcw(K_{s,t})$. 
\end{restatable}
\begin{restatable}[Bernshteyn and Lee~{\cite[Open Problem 68]{BL22}}]{problem}{bla}\label{prob68}
    Is there $n\in \mathbb N$ such that $\tbcw(K_n)>4$? Is it true that in fact $\lim_{n\to\infty} \tbcw(K_n)=\infty$?
\end{restatable}
\begin{restatable}[Bernshteyn and Lee~{\cite[Open Problem 69]{BL22}}]{problem}{blb}\label{prob69}
    Is there a constant $k$ such that $\tbcw(G)\le k$ for all planar graphs $G$? Does $k=4$ work?
\end{restatable}
We answer \zcref{prob68} positively and \zcref{prob69} negatively.
To answer these questions, we prove that every class of graphs of bounded radius-$1$ blind cop-width has bounded treewidth, which can be seen as 
the main contribution of this
paper.

\getkeytheorem{twbcw}

\getkeytheorem{lbcnot}

\getkeytheorem{lbhunt}

For the proof of \zcref{tw}, we want a structure in the graph giving a lower bound for the radius-$1$ blind cop-width.
For that, we introduce the notion of \emph{balanced
minors}.
A minor~$H$ of a graph $G$ is \emph{balanced} if there is a constant $c$ such that every vertex of $H$ corresponds to exactly~$c$ vertices of $G$.
The following proposition, combined with the celebrated Grid Minor Theorem of Robertson and Seymour~\cite{grid}, shows that every graph of sufficiently large treewidth contains a fixed outerplanar graph as a balanced minor.

\getkeytheorem{gridouterplanar}

We then show that graphs containing a large complete binary tree as balanced minors have high blind cop-width, and this achieves the proof of
\zcref{tw} since trees are outerplanar.
Note that having a large complete binary tree as a minor does not guarantee large blind cop-width, as witnessed by \zcref{thm:disprove}.
As a corollary, we obtain a theorem of Bernshteyn and Lee~\cite[Theorem 3]{BL22}, which states that there is a sequence of subcubic trees of arbitrarily large radius-$1$ blind cop-width.
We remark that Althoetmar, Schade, and Schürenberg~\cite{ASS2025} claimed that the complete binary tree of height $d$ has radius-$1$ blind cop-width $(1+o(1))d/2$.
\getkeytheorem{twtree}

We also show that the outerplanarity in \zcref{outerplanar} cannot be relaxed.
\getkeytheorem{iffouterplanar}

If we replace a complete binary tree with a complete graph in \zcref{tw-tree},
then we obtain a more explicit bound as follows.

\getkeytheorem{hbtobcw}
We will show in \zcref{lem:balancedcliqueminor} that if a graph contains $K_{2h-1}$ as a minor, then it contains~$K_h$ as a balanced minor.
This answers \zcref{p66} asymptotically, as $\tbcw(K_n)=\Theta(n)$ and $\tbcw(K_{s,t})=\Theta(\min(s,t))$.

The paper is organized as follows. 
In \zcref{sec:prelim}, we recall some facts
and definitions about the structural graph parameters and Cops and Robber games. 
In \zcref{sec:radius}, we show that blind cop-width parameters with finite radius are functionally equivalent, in the sense that if a class of graphs has bounded radius-$r$ blind cop-width for a finite $r$, then it has a bounded radius-$r'$ blind cop-width for any finite $r'$.
In \zcref{sec:minors}, we introduce the notion of balanced
minors and prove a balanced version of the grid theorem for outerplanar graphs.
In \zcref{sec:bfw}, we study some links between blind cop-width and
blind flip-width and refute \zcref{conj:torunczyk}.
We study the blind
cop-width of subdivisions of graphs in \zcref{subdivisions}. 
In \zcref{connections}, we discuss the links with other games 
studied in the literature.
In \zcref{sec:end}, we discuss some ideas for further work.

\section{Preliminaries}\label{sec:prelim}

For integers $\ell$ and $n$, we write $[\ell,n]$ for the set of integers
$\{\ell, \ell+1,\ldots, n\}$, and $[n]$ for $[1,n]$. All graphs in this paper
are finite, undirected, and simple. For a graph~$G$, we write $V(G)$ for its
vertex set and  $E(G)$ for its edge set. 
For every vertex $v\in V(G)$, let $N_G(v)$ be the set
of all vertices of~$G$ adjacent to~$v$. 
For a set~$A$ of vertices of~$G$, we
define $N_G(A) = \bigcup_{u\in A}N_G(u)\setminus A$ and $N_G[A]=N_G(A)\cup A$. 
A \emph{clique} is a set of pairwise adjacent vertices.
A \emph{complete graph}~$K_n$ is a graph on~$n$ vertices where every two distinct vertices are adjacent. 
A \emph{complete bipartite graph}~$K_{s,t}$ is 
a graph isomorphic to a graph on the vertex set $A\cup B$ 
with $\abs{A}=s$, $\abs{B}=t$, $A\cap B=\emptyset$
such that two vertices are adjacent if and only if one is in~$A$ and the other is in~$B$. 
A graph~$H$ is a \emph{subgraph} of~$G$ if $V(H)\subseteq V(G)$ and $E(H)\subseteq E(G)$, in which case we write $H\subseteq G$. 
It is an \emph{induced subgraph} if $E(H) = E(G)\cap V(H)^2$. 
For any set $X\subseteq V(G)$, $G[X]$ is the subgraph of~$G$ 
\emph{induced by} $X$, and $G\setminus X := G[V(G)\setminus X]$.

\subsection{Minors and subdivisions.}
A \emph{minor model} of a graph~$H$ in a graph~$G$ is a collection
$(X_u)_{u\in V(H)}$ of nonempty pairwise disjoint subsets of~$V(G)$
such that
\begin{itemize}
    \item For every $u\in V(H)$, $G[X_u]$ is connected.
    \item For every edge $\{u,v\}\in E(H)$,  there exist $s\in X_u$ and $t\in X_v$ such that
        $\{s,t\}\in E(G)$.
\end{itemize}
The sets $X_u$ for each $u\in V(H)$ are called the \emph{branch sets} of the minor model.
A graph~$H$ is a \emph{minor} of a graph~$G$ if $G$ contains a minor model of~$H$. 
Alternatively, a graph~$H$ is a minor of a graph~$G$ if $H$ can be obtained from~$G$ by deleting
edges and vertices and contracting edges, where contracting an edge $\{u,v\}$ is an operation that deletes the edge $\{u,v\}$ and identifies~$u$ and~$v$.
It is well known that if $G_1$ is a minor of $G_2$ and $G_2$ is a minor of $G_3$, then $G_1$ is a minor of $G_3$. 

A graph~$H$ is a \emph{subdivision} of a graph~$G$ if it can be obtained from~$G$ by
replacing all edges with internally disjoint paths of nonzero length.
Clearly, if $H$ is a subdivision of a graph~$G$, then $G$ is a minor of $H$.

\subsection{Sparse graph parameters}\label{subsec:sparsity}

A \emph{graph parameter} is a function that maps a graph to a non-negative integer.
A graph parameter~$p$ is extended to classes of graphs by setting $p(\mathscr{C}) = \sup_{G\in
\mathscr{C}}p(G)$. Two graph parameters~$p$ and~$q$ are said to be \emph{functionally equivalent}
if there exist functions $f$, $g$ such that for all graphs~$G$, $p(G) \leq f(q(G))$ and $q(G) \leq g(p(G))$. We write $p\simeq q$ to say that $p$ and $q$ are functionally equivalent.

A class $\mathscr C$ of graphs is \emph{weakly sparse} if there exists a
positive integer~$t$ such that no graph in~$\mathscr C$ contains a subgraph isomorphic to~$K_{t,t}$.

\paragraph{Treewidth.}

A \emph{tree decomposition} of a graph~$G$ is a pair $(T,\mathcal X)$ where $T$ is a tree and $\mathcal X$
is a family of subsets of vertices $X_v\subseteq V(G), v\in V(T)$ satisfying the following three conditions.
\begin{itemize}
    \item $\bigcup_{v\in V(T)}X_v = V(G)$.
    \item For every $e\in E(G)$, there exists $v\in V(T)$ such that $X_v$ contains both ends of $e$.
    \item For every $u\in V(G)$, $T[\{v\in V(T): u\in X_v\}]$ is connected.
\end{itemize}

The \emph{width} of a tree decomposition $(T,\mathcal X)$ is $\max_{v\in V(T)}\abs{X_v}-1$. The \emph{treewidth} of a
graph~$G$, denoted by $\tw(G)$, is the minimum width of a tree decomposition of~$G$. 
The only purpose of having minus $1$ in the definition of width is to make trees have treewidth at most $1$.

\paragraph{Pathwidth.}
A \emph{path decomposition} is a tree decomposition $(T,\mathcal X)$ such that $T$ is a path. The
\emph{pathwidth} of a graph~$G$, denoted by $\pw(G)$, is the minimum width of a path decomposition of~$G$.

\paragraph{Properties.}
Trivially, $\tw(G) \le \pw(G)$.
If $H$ is a minor of~$G$, $\tw(H) \le \tw(G)$ 
and $\pw(H)\le \pw(G)$, see \cite[Proposition 12.3.6]{DiestelGT}.
It is an easy exercise to prove that $\tw(G')=\tw(G)$ for all subdivisions $G'$ of~$G$. 
Taking a subdivision of a graph may increase its
pathwidth by at most two, implied by results on vertex separation numbers and search numbers~\cite{Kinnersley1992,EST1994}.
Classes of graphs of bounded treewidth are weakly sparse, hence so are those of bounded pathwidth.

Binary trees have unbounded pathwidth. 
The \emph{complete binary tree} of height $h$ is a rooted tree such that every non-leaf node has exactly two children and every leaf has distance exactly~$h$ from the root.
\begin{prop}[Scheffler~{\cite[Satz 4 on page 39]{Scheffler1989}}]\label{pathwidth-tree}
	The complete binary tree of height $n$ has pathwidth $\lceil n/2\rceil$.
\end{prop}

\subsection{Cops and robber games}\label{sec:games}

A variety of pursuit-evasion games on graphs can be framed as a number of cops and one
robber standing on vertices of a graph, consecutively moving on it, the cops
trying to catch the robber and the robber trying to escape from the cops. Following the
notation and analysis of Toru\'nczyk~\cite{T23}, we focus on ``helicopter'' games with varying speed
$r\in \mathbb{N}\cup \{\infty\}$ of the robber and write $\copwidth_{r}(G)$ the
minimum number of cops needed to catch a robber of speed~$r$ on a given graph~$G$, using
the following definitions of games and allowed strategies.

\begin{definition}[Helicopter games]
	A \emph{game of radius $r$} on a graph~$G$ is a pair $(\mathcal{A},\mathcal {C})$ of a sequence $\mathcal A$ of vertices
	$a_1, \ldots,a_m \in V(G)$ and a sequence $\mathcal C$ of sets of vertices $C_1, \ldots, C_m \subseteq V(G)$ for some $m\in
	\mathbb{N}$ such that for all $i\in [m-1]$, $a_i\not\in C_i$
    and there exists a path from $a_i$ to~$a_{i+1}$ of length at most $r$ in $G\setminus (C_i\cap C_{i+1})$.

	The game is \emph{winning} (for the cops) if $a_m\in C_m$, and \emph{losing} otherwise. The number
	of cops used is $\max_{i\in [m]}\abs{C_i}$.
\end{definition}

The vertices in~$\mathcal{A}$ represent the positions of the robber moving at speed~$r$
avoiding the cops sitting at positions in~$\mathcal{C}$.
In the visible variant (the usual one), strategies of the cops (defining~$C_{i+1}$ given $C_j$ and $a_j$
for all $j \in [i]$)
are allowed to depend on the position of the robber at the previous turns,
while strategies of the robber are allowed to use the positions of the cops at turns $i$
and $i+1$ to define $a_{i+1}$. The interpretation is that cops move using helicopters
(not being limited to any speed or even to stay in a connected component), and as such
make great noise before landing, allowing the robber to use the information of their
next positions in order to design an escape plan. 

A classical result of Seymour and Thomas~\cite{ST93} states that for all graphs~$G$,
$\copwidth_{\infty}(G) =\tw(G)+1$. Another fundamental parameter of sparsity
theory, degeneracy, can also be linked in a similar way to cops and robber games, as
$\copwidth_{1} = \degeneracy+1$ \cite[Theorem IV.4]{T23}. 

\paragraph{Blind variant.}
In order to characterize pathwidth, 
we introduce the \emph{blind} variant of the game: 
cops' strategies are no longer able to use the position of the robber. 
Formally, $\mathcal{A}$ is now a sequence of sets $A_1,
\ldots, A_m \subseteq V(G)$ such that $A_1=V(G)\setminus C_1$ and
for all $i\in[m-1]$, 
\begin{equation}\label{eq:bcw_r}
A_{i+1} = \{u\in V(G)\setminus C_{i+1}:
\text{$G\setminus (C_i\cap C_{i+1})$ has a path from $A_i$ to $u$ 
of length at most $r$}
\}.
\end{equation} 
The game is winning for the cops if $A_m =\emptyset$.
The \emph{radius-$r$ blind cop-width} of a graph~$G$, denoted by $\bcw_r(G)$, is the minimum number of 
cops so that there is a strategy for the cops to win the blind radius-$r$ cops and robber game.
The sets represent the possible positions of the robber at each turn, and
may also be thought of as sets of vertices infected by a gas spreading at speed~$r$ and
cleaned by the cops. 
We may then say that $V(G)\setminus A_i$ is the set of \emph{clean} vertices.
Since the game is only determined by the strategy of the cops, we might speak of a
\emph{strategy} instead. 
A classical result of Kirousis and Papadimitriou~\cite{KP85}
states that $\bcw_{\infty}(G) = \pw(G) +1$ for every graph~$G$.

\section{Finite radius does not matter much for blind cops}\label{sec:radius}
Before relating blind cop-width with other similar parameters, we prove the following
important lemma that justifies speaking about blind cop-width without specifying the
speed of the robber and has direct implications.

\begin{lem}\label{feq}
	For every  positive integer $r$ and every graph $G$, 
	\[ \bcw_1(G) \le \bcw_r(G)\le 2^{\lceil\log_2 r\rceil} \bcw_1(G)\le 2r \bcw_1(G).\] 
	Consequently, all blind cop-width parameters with finite radius are functionally equivalent.
\end{lem}

\begin{proof}
	Let $G$ be a graph.
    Clearly, $\bcw_r(G) \le \bcw_s(G)$ whenever $r\le s$; if $k$ cops can catch a fast
    robber, they can still catch him if he moves slower.

	We claim that $\bcw_{2r}(G) \le 2 \bcw_r(G)$, from which the statement of
    the lemma follows.
	Here is a strategy for the cops with twice as many cops to catch
    a robber twice as fast: if cops are placed
    at positions $C_i$ at step~$i$ in the speed-$r$ game, the strategy is to place
    cops at positions $C_{2i-1} \cup C_{2i}$ at step $i$ in the speed-$2r$ game.
    Let $R_i$ and $S_i$ be the set of cleaned vertices
    at step $i$ in respectively the $r$ and $2r$ speed games, and let us show that
    $R_{2i} \subseteq S_i$. Initially, $R_0=S_0=\emptyset$. Let $i\ge 1$, and assume that
    $R_{2i-2} \subseteq S_{i-1}$. If $u \in R_{2i}$, either $u\in C_{2i-1}\cup C_{2i}$, in which
    case $u\in S_i$, or already we had $u\in R_{2i-2}\subseteq S_{i-1}$. 
    In the second case, $u$ cannot get infected at step $i$ in the speed-$2r$ game, as it
    would mean there is some vertex $v \notin {S_{i-1}}$, thus not in ${R_{2i-2}}$, with a path
    of length at most $2r$ connecting it to $u$ not hitting $C_{2i-1}\cup C_{2i}$, but then
    $v$ would also infect $u$ in the speed-$r$ game, contradicting the fact that $u\in R_{2i}$.
    Hence $R_{2i}\subseteq S_i$, and the result follows by induction.
\end{proof}

One thing to note is that this property seems to be a consequence of the blindness of the
game, as cops do not need to know what happens between individual rounds to compose
strategies. In particular, this is false for the classical generalized coloring
numbers, also known as non-blind cop-width numbers: 
if~$G$ is a subdivision of~$K_n$ obtained by replacing each edge with a path of length $2r-1$, the cop-width of radius~$r$ of~$G$ will
be less than or equal to~$5$, but cop-width of radius~$2r$ is equal to~$n$.
Since all finite radii are equivalent, 
a class~$\mathscr{C}$ of graphs admits a function $f$ such that the blind cop-width of radius~$r$ of graphs in~$\mathscr C$ is at most $f(r)$ for all $r>0$
if and only if the blind cop-width of radius $1$ of graphs in $\mathscr C$ is bounded. 

\section{Balanced minors and lower bounds}\label{sec:minors}

We introduce balanced minors, which will be useful for obtaining a lower bound for the blind cop-width.
We say that a graph $H$ is a \emph{balanced minor} of a graph~$G$ 
(or $G$ contains $H$ as a balanced minor)
if there is a minor model $(X_u)_{u\in V(H)}$ of $H$ in~$G$ such that all $X_u$ have the same size. 
In this case, we say the minor model is \emph{balanced}.
We write $H\minor G$ if $H$ is a minor of~$G$, and $H \minor_b G$ if $H$ is a balanced minor of~$G$.
It turns out that containing certain balanced minors is sometimes actually equivalent to containing usual
minors, which will be useful in deriving lower bounds for the blind cop-width of families of graphs
defined by usual structural properties.

\subsection{Lower bound lemma}

The only way to prove lower bounds on the blind cop-width of graphs we are aware of is
the following
lemma, of which all proofs of the literature and ours ultimately make use (see
\cite[Proposition 12]{BL22} or \cite[Lemma 1]{XYZZ19}).

\begin{lem}\label{lower_bound}
	Let $G$ be a graph. 
	Let $a$, $k$ be positive integers with $a\le \abs{V(G)}$.
	If 
	\[ 
		\abs{A}+\abs{N_G(A)}\ge a+k
	\] 
	for every set $A$ of vertices with $\abs{A}\ge a$, then $\bcw_1(G)>k$.
\end{lem}
\begin{proof}
	Suppose that $\bcw_1(G)\le k$. 
	Let $(\mathcal{A}, \mathcal{C})$ be a winning game for $k$ cops in $m$ rounds so that $A_1 =V(G)\setminus C_1$ and $A_m = \emptyset$
	and $A_{i+1}=(N_G(A_{i})\cup A_{i})\setminus C_{i+1}$ for all $i\in [m-1]$.
	We take $A_0=V(G)$ and $C_0=\emptyset$ so that $A_{i+1}=(N_G(A_{i})\cup A_{i})\setminus C_{i+1}$ holds for $i=0$ as well. 
	
	Let $i$ be the maximum non-negative integer such that $\abs{A_i}\ge a$.
	Then 
	$a>\abs{A_{i+1}}=\abs{A_i}+\abs{N_G(A_i)}-\abs{C_i}
	\ge \abs{A_i} + \abs{N_G(A_i)}-k $
	and therefore $\abs{A_i}+\abs{N_G(A_i)}< a+k$, 
	contradicting the assumption.
\end{proof}

We give a slight strengthening of this lemma which can make some proofs
easier.
This is also observed by Althoetmar, Schade, and Schürenberg~\cite[Lemma 1]{ASS2025} independently.

\begin{lem}\label{cor:lower_bound}
	Let $G$ be a graph and let $a$, $k$ be positive integers with $a\le \abs{V(G)}$.
	If $\abs{N_G(A)}\ge k$ for all subsets $A$ of $V(G)$ with $\abs{A}=a$, then $\bcw_1(G)>k$. 
\end{lem}

\begin{proof}
    Suppose not. By \zcref{lower_bound},
    there exists a minimal set~$A$ such that $\abs{A}\ge a$ and
    $\abs{A}+\abs{N_G(A)} < a + k$. 
	If $\abs{A}>a$, then for $v\in A$, we have 
	$\abs{A\setminus \{v\}} 
	+\abs{N_G(A\setminus\{v\})}
	\le \abs{A}-1 + \abs{N_G(A)}+1 < a+k$, contradicting the minimality assumption. 
	Therefore, $\abs{A}=a$ and $\abs{N_G(A)}<k$, a contradiction.
\end{proof}

\subsection{Complete graph as a balanced minor}\label{cliques}

Similarly to the \emph{Hadwiger number} $\had(G)=\max\{k\in \mathbb{N}: K_k \minor G\}$, we can
define the \emph{balanced Hadwiger number} $\had_b(G) :=
\max\{k\in \mathbb{N}: K_k \minor_b G\}$.
These two parameters end up being tied to each other.

\begin{lem}\label{lem:balancedcliqueminor}
	For every graph~$G$, 
	$\had(G)/2  \le \had_b(G) \le \had(G) $.
\end{lem}

\begin{proof}
	The second inequality follows from the fact that a balanced minor is itself a minor in the usual sense.
	Thus, it is enough to prove that 
	if $G$ contains $K_{2n-1}$ as a minor, 
	then $G$ contains $K_n$ as a balanced minor.
	Let $(X_i)_{1\le i \le 2n-1}$ be a minor model of~$K_{2n-1}$ in~$G$
	such that $\abs{X_1} \le \abs{X_2} \le \cdots \le \abs{X_{2n-1}}$.
	Let $t=\abs{X_{n}}$.
	For all~$i$ with $\abs{X_i}<t$, 
	let $v_i\in X_{2n-i}\cap N_G(X_i)$.
	Since $G[X_{2n-i}]$ is connected, there exists a subset~$Z_i$ of~$X_{2n-i}$ containing~$v_i$ such that $\abs{Z_i}=t-\abs{X_i}$
	and $G[Z_i]$ is connected.
	For all $j$ with $\abs{X_j}=t$, we define $Z_j=\emptyset$.
	Now, for all $i$ with $\abs{X_i}\le t$, we define $Y_i=X_i\cup Z_i$.  
	Then $(Y_i)_{1 \le i \le n}$
	is a balanced minor model of $K_{n}$.
\end{proof}

The bound in \zcref{lem:balancedcliqueminor} is tight for all odd $n:=\had(G)\ge 5$. Indeed, for $n\ge 4$, let $G_n$ be the subdivision
of $K_{n}$ obtained by subdividing each edge 
in a fixed $K_{\lceil n/2\rceil}$  subgraph at least $2n$ times. It can be checked that $\had_b(G) =
\lfloor n/2\rfloor+1$,
while $\had(G)=n$.
Here we need $n\ge 4$ so that each branch set of a clique minor model of order greater than $\lfloor 
n/2\rfloor+1$ is forced to contain a vertex of degree at least $3$.

The following proposition will allow us to give a linear lower bound to the blind
cop-width by the size of a largest clique minor. 
We introduce the following notations that will be reused later. For  a minor model $(X_i)_{i\in[k]}$ and a set $A$ of vertices, 
\begin{itemize}
    \item 
	let $F_A=\{i\in[k]: X_i\cap A = X_i\}$
	be the set of indices of \emph{full} branch sets with respect to~$A$, 
    \item let $E_A=\{i\in[k]:  X_i\cap A = \emptyset\}$ be the set of indices of \emph{empty} branch sets with respect to~$A$, and 
    \item let $M_A=\{i\in[k]:  X_i\cap A\neq\emptyset,~ X_i\cap A\neq X_i\}$ be the set of indices of \emph{mixed} branch sets with respect to~$A$.
\end{itemize} 
In particular, $|F_A|+|E_A|+|M_A| = k$.

\begin{prop}[store=hbtobcw]\label{prop:hbtobcw}
    If a graph~$G$ contains $K_h$ as a balanced minor, then 
    \[ \bcw_1(G)\ge \min(h, (h+3)/2).\]
\end{prop}
\begin{proof}
    Let $(X_i)_{i\in [h]}$ be a balanced minor model of~$K_h$ in~$G$. 
	Let $s=\abs{X_1}=\abs{X_2}=\cdots=\abs{X_h}$.
	We may assume that $\bigcup_{i=1}^h X_i = V(G)$ by deleting unnecessary vertices, because the blind cop-width clearly does not increase when taking subgraphs.
    If $h\le 2$, then $\bcw_1(G)\ge h$ trivially. 
	If $s=1$, then $G=K_h$ and $\bcw_1(G)=h$. 
	Thus we may assume that $h\ge 3$ and $s\ge2$.

	Let $k=\lfloor h/2\rfloor$ and 
	let $A$ be a set of vertices such that $\abs{A}= (s-1)k+1$. 
    As $\lceil (h+3)/2\rceil=k+2$, 
	by \zcref{cor:lower_bound}, it is enough to prove that 
	\begin{equation}\label{eq:a+na}
		\abs{N_G(A)}\ge k+1 .
	\end{equation}
	Let $F_A$, $E_A$, $M_A$ be the set of indices of full, empty, and mixed branch sets in $(X_i)_{i\in [h]}$ with respect to $A$.

	Suppose that $F_A=\emptyset$. Then $X_i$ is mixed whenever $X_i\cap A\neq\emptyset$. Thus, %
	\(  
		\abs{M_A}\ge \frac{1}{s-1}\abs{A} > k
	\).
	Since each mixed branch set induces a connected subgraph and partially intersects $A$, $A$ has a neighbor in each mixed branch set, thus $\abs{N_G(A)} \ge \abs{M_A} \ge  k+1$, implying \eqref{eq:a+na}.

	Thus we may assume that $F_A\neq\emptyset$.
	Since a full branch set is connected to every empty branch set, we have $\abs{N_G(A)}\ge \abs{M_A}+\abs{E_A} = h - \abs{F_A}$.
	Since $\abs{F_A}\le \frac{1}{s}\abs{A} \le k - \frac{k-1}{s}$, 
	we have $\abs{N_G(A)}\ge h-k+\frac{k-1}{s}$.
	If $h\ge4$, then $k-1>0$ and so $\abs{N_G(A)}\ge h-k+1 \ge k+1$.
	If $h=3$, then $k=1$ and $\abs{N_G(A)}\ge h -k = k+1$.
	Thus in both cases, we have \eqref{eq:a+na}.
	This completes the proof.
\end{proof}

\begin{figure}
	\begin{center}
		\begin{tikzpicture}
			\tikzstyle{v}=[circle, draw=black, solid, fill=black, inner sep=0pt, minimum width=2pt]
			\tikzstyle{vv}=[circle, draw=red, solid, fill=red, inner sep=0pt, minimum width=1pt]

			\foreach \t in {1,2,3,4,5,6,7,8} 
			{ 
				\node[v] (v\t) at (\t*360/8-0.5*360/8+90:4) {};
			}
			\foreach \i in {1,2,3,4}{
				\pgfmathtruncatemacro{\nexti}{\i+1}
				\foreach \j in {\nexti,...,8} {
					\draw (v\i)--(v\j);
				}
			}
			\foreach \i in {1,2,3,4}{
				\draw[red,ultra thick] (v\i) 
				--
				+(-4,0)node [v]{};
			}
			\foreach \i in {5,6,7}{
				\pgfmathtruncatemacro{\nexti}{\i+1}
				\foreach \j in {\nexti,...,8} {
					\draw [red,ultra thick] (v\i)--(v\j);			
				}
			}
		\end{tikzpicture}
	\end{center}
	\caption{A graph having radius-$1$ blind cop-width $6$ or $7$ and containing $K_8$ as a balanced minor. Thick red edges have been subdivided $2^{12}$ times.}	\label{fig:g8}

\end{figure}
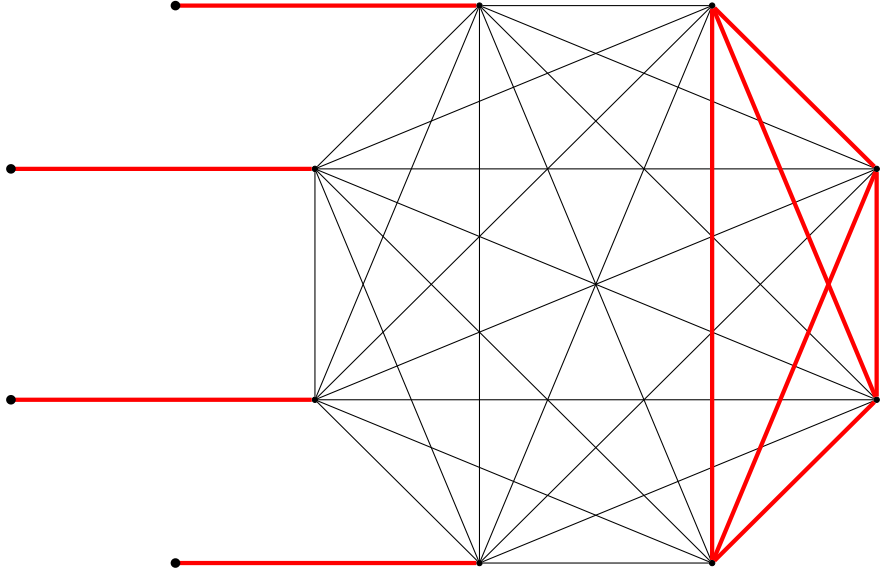

The bound in \zcref{prop:hbtobcw} is also close to optimal.
Let $G$ be a graph obtained from~$K_{2t}$ on the vertex set $[2t]$ by subdividing all edges $\{i,j\}$ with $t+1\leq i<j\leq 2t$ at least $2^{t(t-1)}$ times, and appending a sufficiently long path to each vertex $i\in [t]$. See \zcref{fig:g8} for an illustration.
It can be checked that $h(G)=h_b(G)=2t$
and therefore $\bcw_1(G)\ge t+2$ by \zcref{prop:hbtobcw}.
On the other hand, the radius-$1$ blind cop-width is at most $t+3$.
The strategy for the $t+3$ cops is as follows.
\begin{enumerate}[label=\rm(\arabic*)]
	\item At first, $t$ cops camp the branch vertices in $[t]$ until step (4).
	\item Two other cops clean the long paths appended to the branch vertices in $[t]$.
	\item The $3$ cops can perform a cleaning process so that in the end, 
a ball of radius $1$ around each of the branch vertices in $[2t]\setminus [t]$
is clean. Let $(x_1,y_1)$, $(x_2,y_2)$, $\ldots$, $(x_{t(t-1)},y_{t(t-1)})$ be the list of all pairs $(x,y)$ of distinct elements of $[2t]\setminus [t]$
where $x_1\le x_2\le \cdots \le x_{t(t-1)}$.
	\begin{itemize}
		\item From $i=1$ to $i=t(t-1)$, one cop occupies $x_i$ 
			and two other cops clean $2^{t(t-1)-i}$ internal vertices closest to $x_i$ of the path of $G$ corresponding to the edge $x_iy_i$ of $K_{2t}$.
			This takes time at most $2^{t(t-1)-i}$.
	\end{itemize}
	Note that at the end of this process, at each branch vertex  in $[2t]\setminus [t]$, a ball of radius $1$ around it is clean.
\item The $t$ cops move from $[t]$ to $[2t]\setminus [t]$. By the previous step, the branch vertices in $[t]$ remain clean.
\item The $2$ free cops clean the remaining degree-$2$ vertices.
\end{enumerate}

\subsection{Highly connected minors with many cliques}

The bounds on the blind cop-width of complete
balanced minors are not
sufficient to prove our main theorem in the general case, but
can readily be used in certain highly-connected settings, as
we shall see in this subsection. In particular, when 
applicable, those bounds are much better than the ones
purely in terms of treewidth.

\begin{prop}\label{lem:disjointcliques}
	If a graph~$G$ has a $k$-connected minor~$H$ having 
	$k$ pairwise disjoint cliques of $k$ vertices,
	then $G$ contains $K_k$ as a balanced minor.
\end{prop}
\begin{proof}
	Let $C_1$, $C_2$, $\ldots$, $C_{k}$ be pairwise disjoint cliques in~$H$ such that $\abs{C_1}=\cdots=\abs{C_{k}}=k$.
	Let $(X_i)_{i\in V(H)}$ be a minor model of~$H$ in $G$.
	We may assume that $m=\max_{j\in C_1}\abs{X_j}\le \max_{j\in C_i}\abs{X_j}$ for all $i\in [k]$. 
	Let $S=\{i\in C_1: \abs{X_i}<m\}$ and 
	let $T=\{i\in \bigcup_{j=2}^{k} C_j: \abs{X_i}\ge m\}$.
	By the assumption on $C_1$, $\abs{S}\le k-1$ and $\abs{T}\ge k-1$.
	Since $H$ is $k$-connected, 
	by Menger's theorem,
	$H\setminus (C_1\setminus S)$ has $\abs{S}$ pairwise vertex-disjoint paths from $S$ to $T$.
	Each of those paths corresponds to a path in~$G$ joining~$X_i$ for some $i\in S$ to $X_j$ for some $j\in T$.
	Since each $X_j$ for $j\in T$ has at least $m$ vertices,
	we can use these paths to increase the size of $X_i$ for every $i\in S$ 
	to have exactly $m$ vertices. 
	Thus $G$ contains $K_k$ as a balanced minor.
\end{proof}

A \emph{$k$-tree} is a graph built from $K_{k+1}$ by repeatedly adding a new vertex adjacent to all vertices in a clique of size $k$. 
We show that every $k$-tree with many pairwise disjoint cliques of size $k$ has topological blind cop-width at least $(k+3)/2$.

\begin{prop}\label{prop:ktree}
	Let $k\ge 3$ be an integer. 
	Let $G$ be a $k$-tree with $k$ pairwise disjoint cliques of size~$k$.
	Then every subdivision of $G$ has radius-$1$ blind cop-width at least $(k+3)/2$. 
\end{prop}
\begin{proof}
	Since every $k$-tree is $k$-connected, by
    \zcref{lem:disjointcliques}, 
	$G$ contains $K_k$ as a balanced minor. 
	By \zcref{prop:hbtobcw}, the radius-$1$ blind cop-width of $G$ is at least $(k+3)/2$.
\end{proof}

Here is an easy corollary.

\begin{cor}
	For every integer $k\ge 3$, 
	every subdivision of the $k$-th power of a path on $k^2$ vertices has radius-$1$ blind cop-width at least $(k+3)/2$. 
	\qed
\end{cor}

\subsection{Outerplanar graph as a balanced minor}

Applying similar vertex exchange techniques to large grid minors 
instead of large clique minors, 
we can still find a useful type of balanced minors: outerplanar balanced minors. 
These outerplanar minors are particularly valuable because
they can force the blind cop-width to be arbitrarily large, as we show in the next subsection. 

\begin{definition}
	The $n\times n$ grid is the graph whose
	vertex set is ${[n]}^2$ where $(k,\ell)$ and $(k',\ell')$ are adjacent if and
	only if $\abs{k-k'}+\abs{\ell-\ell'}=1$.
\end{definition}

Robertson and Seymour~{\cite[(1.3)]{RST1994}} proved that 
every $n$-vertex planar Hamiltonian graph is a minor of an $n\times n$ grid.
As observed by Gavoille and Hilaire~\cite{GH2023}, the proof of Robertson and Seymour implies that every $n$-vertex outerplanar graph~$G$ is a minor of the upper half of the $n\times n$ grid, that is the subgraph of $n\times n$ grid induced on $\{(i,j): 1\le i\le j\le n\}$, where $n$ vertices on the diagonal are mapped to the vertices of~$G$. For completeness of this paper, we include a proof.

\begin{lem}\label{lem:embed-outerplanar}
	For every $n$-vertex outerplanar graph~$G$, 
	the $n\times n$ grid has a minor model $(X_v)_{v\in V(G)}$ of~$G$ such that 
	\begin{enumerate}[label=\rm(\alph*)]
		\item for every vertex~$v$ of~$G$, there is a unique $k\in [n]$ with $(k,k)\in X_v$ 
		and 
		\item for every $k\in [n]$, 
		$X_k$ contains no vertex~$(a,b)$ with $a>b$.		
	\end{enumerate}
\end{lem}
\begin{proof}
	We follow the proof of Robertson and Seymour~{\cite[(1.3)]{RST1994}}.
	Let us fix an embedding of~$G$ on the plane, where all vertices are on the infinite face.
	We may assume that $G$ is maximally outerplanar, so all inner faces are triangles. 
	Let $v_1$, $v_2$, $\ldots$, $v_n$ be the vertices in the cyclic order they appear on the boundary of the infinite face. Since $G$ is maximally outerplanar, $v_1v_2\cdots v_n v_1$ is a Hamiltonian cycle of~$G$.

	For $k\in [2,n]$, let $\ell(k)\in [k-1]$ be the minimum index such that $v_{\ell(k)}$ is adjacent to $v_k$.
	For $k\in [n-1]$, let $r(k)\in [k+1,n]$ be the maximum index such that $v_{r(k)}$ is adjacent to $v_k$. 
	For $k\in [2,n-1]$, let 
	\[ 
		X_k=\{ (i,k): \ell(k)<i \le k\} \cup \{ (k,j): k\le j<r(k)\}.
	\] 
	Let 
	$X_1=\{ (1,j):  1\le j< r(1)\}$ and 
	$X_n=\{ (i,n):  \ell(n)\le i \le n\}$.
	Then $G[X_k]$ is connected for all $k\in [n]$ and $X_1$, $X_2$, $\ldots$, $X_n$ are pairwise disjoint.
	Since all inner faces of $G$ are triangles, if $k<k'$, $r(k)=k'$, and $\ell(k')=k$, then $k=1$ and $k'=g$. 
	Now we deduce that if $v_k$ and $v_{k'}$ are adjacent in~$G$ with $k<k'$, 
	then 
	$(k+1,k')\in X_{k'}$, $(k,k'-1)\in X_k$, and $(k,k')$ belongs to either~$X_k$ or~$X_{k'}$, implying that the $n\times n$ grid has an edge between $X_k$ and $X_{k'}$. 
	Therefore $(X_k)_{k\in [n]}$ is a minor model of $G$ in the $n\times n$ grid. 
\end{proof}

\begin{prop}[store=gridouterplanar]\label{outerplanar}
	If a graph contains the $n^2\times 4n$ grid as a minor, 
	then it contains every $n$-vertex outerplanar graph as a balanced minor. 
\end{prop}

\begin{figure}
    \centering

	\begin{tikzpicture}[scale=0.8]
		\tikzstyle{v}=[circle, draw=black, solid, fill=black, inner sep=0pt, minimum width=2pt]
		\draw [fill=gray!30] (0,2) rectangle++(16,2);
		\foreach \y in {0,...,4} {
			\draw (0,2*\y)--+(2*8,0);

		}
		\foreach \x in {0,...,7}{
			\draw(2*\x,0)--+(0,8);
			\draw(4+.25*\x,6)--+(0,2);
			\draw(4,6+.25*\x)--+(2,0);
			\node[v] (v\x) at (4+.25*\x, 6.25+.25*\x){};
		}
		\draw (16,0) --+(0,8);
		\foreach \x/\y in {0/1,2/3,4/5,6/7,0/3,0/7,4/7,2/7}{
			\node[v] (v\x\y) at ((4+.25*\x, 6.25+.25*\y){};
		}
		\draw [ultra thick](v01)--(v03)--(v07)--(v27)--(v67);
		\draw [ultra thick](v7)--(v67)--(v6);
		\draw [ultra thick](v2)--(v23)--(v3);
		\draw [ultra thick](v4)--(v45)--(v5);
		\draw [ultra thick](v0)--(v01)--(v1);
		\draw [ultra thick] (v03)--(v23);
		\draw [ultra thick] (v47)--(v45);
		\foreach \x in {0,1,2}{
			\draw [thick] (v\x)--(4+.25*\x,6-.25*\x) node[v]{} --(2*\x+.25,6-.25*\x) node[v]{}--(2*\x+.25,4);
		}
		\foreach \x in {3,...,7}{
			\draw [thick] (v\x)--(4+.25*\x,4.25+.25*\x) node[v]{} --(2*\x+.25,4.25+.25*\x) node[v]{}--(2*\x+.25,4);
		}	
		
	\end{tikzpicture}
    \caption{A binary tree as a balanced minor of a grid minor model.}
    \label{fig:tree-in-grid}
\end{figure}
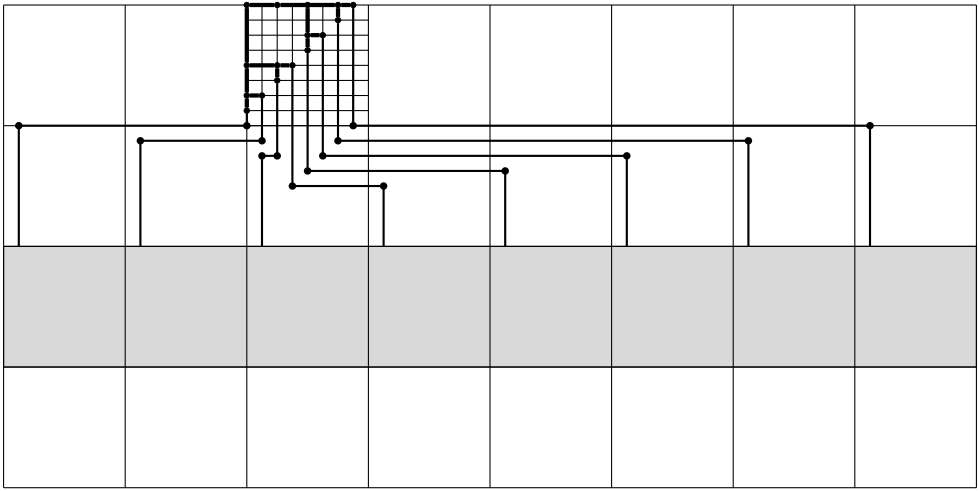

\begin{proof}
	Let $G$ be an $n$-vertex outerplanar graph.
	We may assume that $n\ge 4$, because it is easy to check for $n<4$.
	By~\zcref{lem:embed-outerplanar}, $G$ is a minor of the upper half of the $n\times n$ grid.

    Now let $H$ be a graph containing a minor model $(X_{i,j})_{i,j\in [n^2]}$ of the $n^2\times 4n$ grid.
    Let $I_a := [(a-1)n+1, an]$ for $a\in
	[n]$.
	We then consider the minor model of the $n\times n$ grid in $H$
	 where each branch set is the union of the branch sets of some $n\times n$ subgrid of the $n^2\times n^2$-grid model.  %
	Let us choose $a\in[n]$, $b\in[4]$ such that 
	$\sum_{i\in I_a, j\in I_b} \abs{X_{i,j}}$ is minimized, and let $H'$.
    
	Either
	$b\ge 3$ or $b\le 2$, 
	so one side of~$H'$ faces at least two other rows.
	Let us assume by symmetry that $b\ge 3$ so that it is the bottom side. 
	Let  $(Y_i)_{i\in [n]}$ be the minor model of~$G$ in
 	the upper half of the $n\times n$ grid 
	given by \zcref{lem:embed-outerplanar} 
	where 
	the branch set $Y_i$ intersects the diagonal vertices exactly on the vertex $(i,i)$. 
	Then we can find a minor model $(Z_i)_{i\in [n]}$ of $G$ in $H'$ with~$Z_i$ intersecting the bottom row exactly at $X_{i,1}$, by
	extending each set with a vertical path from the bottom row to the diagonal vertex.
    It is then easily seen that one
	can link the $i$-th vertex to the cell spanned by $I_i\times I_{b-2}$ using pairwise vertex-disjoint
	paths (see \zcref{fig:tree-in-grid}). This allows us to define new branch sets $\Omega_i$ for each vertex $v_i$ as the union of $Z_i$ and a connected subgraph of the path and
 the cell spanned by $I_i\times I_{b-2}$, and still have a minor model of~$G$. Since the cells $I_i\times I_{b-2}$ have more vertices than $H'$ by minimality of $H'$, we can choose the connected subgraphs in a way that ensures that all $\Omega_i$'s have the same size, so that we get a balanced minor model of~$G$ in $H$. %
\end{proof}

This result is in fact optimal in the sense that outerplanar graphs are exactly the
graphs which, when forbidden as balanced minors, also forbid arbitrarily large grid minors.

\begin{prop}[store=iffouterplanar]
	For every non-outerplanar graph~$G$ and every integer~$k$, there exists a graph that contains the $k\times k$ grid as a minor 
    and does not contain $G$ as a balanced minor.
\end{prop}
\begin{proof}
	Let $G$ be a non-outerplanar graph,
	and for a positive integer $n$, let $H_n$ be a subdivision of the $n\times n$ grid where 
	all horizontal edges of the $i$-th row 
	(of the form $\{(b,i),(b+1,i)\}$) are
	subdivided $(2n\abs{V(G)})^i-1$ times. 
    We embed it in the plane by placing the original
	vertices at corresponding integer coordinates and the ones obtained from subdividing 
	evenly in-between, with edges drawn as straight line segments connecting the vertices.

	Assume by contradiction that $G$ is a balanced minor 
	of~$H_n$, witnessed by the minor model $(X_u)_{u\in V(G)}$. 
    By continuously contracting each $X_u$ to a point in the embedding of $H_n$, we obtain an embedding of $G$.
	Since $G$ is not outerplanar, not all vertices of $G$ are on the boundary of the outer face in the embedding of~$G$. So there must be a vertex $u\in
	V(G)$ such that $u$ is separated from the outer face by a cycle $C$ in $G$. 
	Then $H_n[\bigcup_{v\in V(C)}X_v]$ contains a horizontal edge at some
	row of index~$i$ strictly greater than the index of any row intersected by $X_u$,
	thus $\sum_{v\in V(C)}\abs{X_v} \ge (2n\abs{V(G)})^i$. 
	By the pigeonhole principle, 
	there exists some vertex $v\in
	V(C)$ with $\abs{X_v} \ge (2n)^i\abs{V(G)}^{i-1}$. But since $X_u$ only intersects rows with indices in $[i-1]$, 
	\begin{align*}
		\abs{X_u} &\le \sum_{j=1}^{i-1}n(2n\abs{V(G)})^j \\ &\le
		n\abs{V(G)}^{i-1}\sum_{j=1}^{i-1}(2n)^j \\
		&\le \abs{V(G)}^{i-1} \frac{n}{2n-1}((2n)^i-1) \\ &<
		\abs{V(G)}^{i-1}(2n)^i.
	\end{align*}
 
	Hence $\abs{X_u}<\abs{X_v}$ and the minor model is not balanced, contradicting the assumption.
\end{proof}

\subsection{Complete binary tree as a balanced minor}\label{main_thm}

We aim to prove that graphs of bounded radius-$1$ blind cop-width have bounded treewidth.
We do this by showing the existence of certain outerplanar balanced minors is sufficient to give a lower bound on the radius-$1$ blind cop-width.
Specifically, we show this for complete binary tree balanced minors (contrary to the case of complete binary tree minors in general, which do not meaningfully constrain the blind cop-width for any radius).

We begin by considering the subgraph induced by the union of the branchsets of a balanced minor model of a large complete binary tree. Assuming for contradiction that the radius-$r$ blind cop-width is small, we again use \zcref{cor:lower_bound} to get some set of vertices $A$ with desired cardinality (a specific multiple $p\cdot t$ of the size $t$ of the branchsets, to be calculated later) and small
neighborhood with respect to the blind cop-width. We label each branchset of the balanced minor as either full, empty, 
or mixed with respect to $A$, as per \zcref{cliques}.
The first observation is that if there are many mixed branchsets then the neighborhood of $A$ will be large, giving the desired contradiction.
Thus we may assume that almost all branchsets are full or empty, allowing us to approximate the number of full and empty branchsets up to a small additive term (in terms of the ratio $p$ which we selected).
By choosing $p$ to be roughly a third of a sufficiently large power of $4$, we can then show that many empty branchsets are adjacent to full branchsets. This also implies that $A$ has large neighborhood, yielding the desired contradiction.

\begin{prop}[store=twtree]\label{tw-tree}
	There is a function $g:\mathbb{N}\to \mathbb{N}$ such that for all $k\in\mathbb{N}$,
    every graph having a balanced minor model of a 
	complete binary tree of at least $g(k)$ nodes has radius-$1$ blind cop-width larger than~$k$.
\end{prop}

The following lemma is an easy observation on the number of vertices in a subtree based on the location of neighbors and roots.
\begin{lem}\label{lem:subtree-counting}
	Let $T$ be a complete rooted binary tree and let $T'$ be a subtree of $T$. 
	For a node~$v$ of $T$, let $h(v)$ be the minimum distance from $v$ to a leaf of $T$.
	Let $u$ be the unique vertex of~$T'$ with maximum $h(u)$.
	Let $N$ be the set of vertices $w\in N_T(V(T'))$ with $h(w)<h(u)$. 
	Then 
	\[ 
		\abs{V(T')}=2^{h(u)+1}-1 -\sum_{w\in N} (2^{h(w)+1}-1).
	\]
    \qed
\end{lem}

Let us recall the non-adjacent form of a positive integer, mostly from~\cite[page 98]{HMV2004}.
A \emph{signed digit representation} of an integer $k$ is a sum of the form $\sum_{i=0}^{\ell-1} k_i 2^i$, where $k_i\in\{0,1,-1\}$. 
A \emph{non-adjacent form (NAF)} of a non-zero integer $k$ is 
an expression $k=\sum_{i=0}^{\ell-1} k_i2^i$ where $k_i\in\{0,1,-1\}$, $k_{\ell-1}\neq0$, and $k_{i}k_{i-1}=0$ for all $i\in [\ell-1]$. Reitwiesner~\cite{Reitwiesner1960} proved that every non-zero integer~$k$ has a unique non-adjacent form
and furthermore the non-adjacent form of $k$ minimizes 
the number of non-zero digits in a signed digit representation of~$k$. Let us write this number $w(k)$ and call it the \emph{NAF weight} of~$k$. Assume $w(0)=0$.
This property implies that $w(x+y)\ge w(x)-w(y)$ for integers $x$ and $y$, because $w(x+y)+w(y)=w(x+y)+w(-y)\ge w(x+y -y) = w(x)$.
\begin{lem}\label{lem:naf-lemma}
	There is a function $g:\mathbb N\to\mathbb N$ such that 
	for every positive integer $k$, 
	there is no integer in $[g(k)-k, g(k)+k]$ that can be written as a sum of at most $3k$ integers, each of the form $\pm(2^\ell-1)$ for some positive integers~$\ell$.
\end{lem}
\begin{proof}
	Let $S_k$ be the set of all integers that can be written as a sum of at most $3k$ integers of the form~$\pm(2^\ell-1)$ for some positive integers~$\ell$.

	Suppose that $x\in S_k$. Then $x=\sum_{i=1}^m \varepsilon_i (2^{\ell_i}-1)$ with $m\le 3k$ and $\varepsilon_i\in \{1,-1\}$ for all $i\in[m]$. 
	Let $C=\sum_{i=1}^m \varepsilon_i$. 
	Let $w(n)$ denote the NAF weight of $n$. 
	Then $w(x+C)\le m\le 3k$ and $\abs{C}\le 3k$.
	
	Thus, it is enough to prove that for all integers $j$ and $C$ with $\abs{j}\le k$ and $\abs{C}\le 3k$, we have $w(g(k)+j+C)> 3k$. By taking $C'=j+C$, it suffices to show that for all integers $C'$ with $\abs{C'}\le 4k$, we have $w(g(k)+C')>3k$.
	By the inequality on $w$, 
	\[ w(g(k)+C')\ge w(g(k))-w(C').\]
	Let $W(k)=\max \{w(C'): \abs{C'}\le 4k\}$.
	Note that the minimum positive integer $A_w$ with NAF weight $w$ is $2^{2w-2}-2^{2w-4}-2^{2w-4}-\cdots-2^2-2^0=\frac{1}{3}(2^{2w-1}+1)$.
	So if $\frac{1}{3}(2^{2w-1}+1)\le 4k$, then $W(k)\ge w$.
	Thus, $W(k)= \lfloor \frac12(\log_2(12k-1)+1)\rfloor$.

	So, it is enough to make sure that $w(g(k))\ge 3k+W(k)+1$.
	We choose $g(k)=2^0+2^2+2^4+\cdots+2^{2(3k+W(k))} = \frac{1}{3}(4^{3k+W(k)+1}-1)$.
	Then by the construction, $w(g(k))=3k+W(k)+1$, concluding the proof.
\end{proof}
\zcref{lem:naf-lemma} uses $g(k)=\Theta(k\cdot64^k)$, which we did not try
to optimize.
\begin{proof}[Proof of \zcref{tw-tree}]
	Let $g:\mathbb N\to \mathbb N$ be the function given by \zcref{lem:naf-lemma}.
    Let $k$ be a positive integer. 
    Let $G$ be a graph having a balanced minor model~$(X_v)_{v\in V(T)}$ of a complete rooted binary tree~$T$ with $|V(T)|\geq g(k)$.
	Let $g=g(k)$ and let $t=\abs{X_v}$ for all $v\in V(T)$.
    Note that by taking an induced subgraph of $G$, we only make it harder for the robber to escape, without diminishing the capabilities of the cops, so we may assume that $V(G)=\bigcup\{X_v:v\in V(T)\}$.
    
    By \zcref{cor:lower_bound}, it is enough to show that 
	$\abs{N_G(A)}\ge k$
    for every subset $A$ of $V(G)$ with $\abs{A}= gt$.
    Fix such a set $A$, and let $F_A=\{ v\in V(T): X_v\subseteq A\}$, $E_A=\{v\in V(T): X_v\cap A=\emptyset\}$, and $M_A=V(T)\setminus (F_A\cup E_A)$.
	Since each $G[X_v]$ is connected, for each $v\in M_A$, $X_v$ contains a vertex in $N_G(A)$ and therefore $\abs{N_G(A)}\ge \abs{M_A}$. 
	Therefore we may assume that $\abs{M_A}<k$. 

	Then $\abs{A}\le \abs{F_A}t + \abs{M_A}t$ and therefore 
	$\abs{F_A}\ge \frac{1}{t}\abs{A}-\abs{M_A} \ge g-k+1$. 
    Moreover $\abs{F_A}\le \frac{1}{t}\abs{A} = g$.

	We claim that for each node $v\in N_T(F_A)$, $X_v$ contains at least one vertex in $N_G(A)$. If $v\in E_A$, then $X_v$ contains a neighbor of some vertex of $X_w$ with $w\in F_A$ by the definition of a minor model.
	If $v\in M_A$, then since $G[X_v]$ is connected, $X_v$ contains a vertex in $N_G(A)$. 
	Therefore $\abs{N_G(A)}\ge \abs{N_T(F_A)}$.

	Let $T_1$, $T_2$, $\ldots$, $T_m$ be the connected components of $T[F_A]$.
	By \zcref{lem:subtree-counting}, $\abs{V(T_i)}$ can be written as a sum of at most $\abs{N_T(V(T_i))}$ integers of the form $\pm(2^{\ell}-1)$ for a positive integer~$\ell$, except for possibly one of the subtrees $T_i$ containing the root of~$T$, that requires $\abs{N_T(V(T_i))}+1$ integers.
	Note that $\sum_{i=1}^m \abs{N_T(V(T_i))} \le 3\abs{N_T(F_A)}$ because a node of~$T$ can be a neighbor of at most three other nodes in~$T$.
	Thus, we deduce that 
	$\abs{F_A}$ can be written as a sum of at most $3\abs{N_T(F_A)}+1$ integers of the form $\pm(2^{\ell}-1)$ for a positive integer~$\ell$.

	As $g-k<\abs{F_A}\leq g$, 
	by \zcref{lem:naf-lemma}, 
	$\abs{F_A}$ cannot be written as a sum of at most $3k$ integers of the form $\pm(2^{\ell}-1)$ for a positive integer~$\ell$.
	Therefore $3\abs{N_T(F_A)}+1\ge 3k+1$ and therefore $\abs{N_G(A)}\ge \abs{N_T(F_A)}\ge k$.
\end{proof}

Since trees are outerplanar, they appear as a balanced minor in every graph having a sufficiently large grid minor by
\zcref{outerplanar}, and as such in every graph of large treewidth by the
famous grid minor theorem.

\begin{thm}[Robertson and Seymour~\cite{grid}]\label{grid}
	For every integer $k$ there is an integer $f(k)$ such that every graph of treewidth
	at least $f(k)$ has a $k\times k$ grid minor.
\end{thm}

\begin{thm}[store=twbcw]\label{tw}
	For every integer $k$, there is an integer $g(k)$ such that every graph of treewidth
	at least $g(k)$ has radius-$1$ blind cop-width larger than~$k$.
\end{thm}
\begin{proof}
	Let $f_1$ be obtained from \zcref{tw-tree}, $f_2(k)=k^2$ from 
	\zcref{outerplanar}, and $f_3$ from \zcref{grid}. 
	Then every graph with
	treewidth at least $f_3(f_2(f_1(k)))$ has radius-$1$ blind cop-width larger than~$k$.
\end{proof}
Chuzhoy and Tan~\cite{CT2020} showed that the function $f_3$ from \zcref{grid} can be taken as $O(k^{9} \log^{O(1)}(k))$. As the function $f_2$ from \zcref{outerplanar} is a quadratic polynomial and our function $f_1$ from \zcref{tw-tree} is a single exponential function, we can deduce that the function $g$ in \zcref{tw} can be taken as a single exponential function.

\section{Blind flip-width}\label{sec:bfw}

Motivated by a conjecture on first-order model checking, 
Toru\'nczyk~\cite{T23} introduced the flip-width, aiming to be a graph parameter suited for dense graphs defined in terms of games.
For an integer~$k$, a \emph{$k$-flip} of a graph~$G$ is defined by a pair $(\mathcal{P},F)$ of a partition~$\mathcal{P}$ of~$V(G)$ into at most~$k$ parts and a graph~$F$ with vertex set~$\mathcal{P}$, 
where $F$ is allowed to have loops. 
The graph obtained after \emph{flipping~$G$ according to
$(\mathcal{P},F)$} is then the graph~$H$ with vertex set~$V(G)$ such that an edge
is present between two distinct vertices~$u$ and~$v$ if and only if for their
respective parts~$P$ and~$Q$ in~$\mathcal{P}$, either $\{u,v\}\in E(G)$ and
$\{P,Q\}\not\in E(F)$ or $\{u,v\}\not\in E(G)$ and $\{P,Q\}\in E(F)$. 
The \emph{flip-width game} is played on a graph~$G$ between a robber and a flipper.
Let $G_0=G$.
At the $i$-th round, the flipper announces a flip that will be performed to~$G$, 
the robber moves at speed at most~$r$ in~$G_{i-1}$,
and then the flipper performs the announced flip to~$G$ to obtain $G_i$. 
The robber loses if he is trapped on an isolated vertex in $G_i$ for some $i\ge 0$. 
The \emph{flip-width of radius~$r$}, denoted by $\fw_r$, 
is the minimum~$k$ such that there is a winning game for the flipper using $k$-flips.
Here, the number of parts of a partition plays the same role as the number of
cops for cop-width. 
A class~$\mathscr{C}$ of graphs is \emph{of bounded flip-width}
if for every finite radius $r\in \mathbb{N}$, there is some $c_r\in \mathbb{N}$
such that $\fw_r(G) \leq c_r$ for all $G\in \mathscr{C}$. Similarly to cop-width, one
can then define the \emph{blind flip-width of radius~$r$} denoted by
$\bfw_r$ by further restricting the flipper to strategies oblivious to the
position of the robber, and define the notion of classes \emph{of bounded blind
flip-width} as classes~$\mathscr{C}$ admitting a function $f:\mathbb {N}\to\mathbb{N}$ such that $\bfw_r(G)\le f(r)$ for all $r\in \mathbb{N}$ and $G\in \mathscr{C}$.
For $r=\infty$, the parameters $\fw_\infty$ and $\bfw_\infty$ 
are functionally equivalent to the well-studied parameters \emph{clique-width} and
\emph{linear clique-width} \cite[Theorems V.17 and XI.9]{T23}, respectively. 

\subsection{Blind flip-width of weakly sparse graphs}

The motivating question of this work was to answer the following conjecture of Toru\'nczyk.
\getkeytheorem{torunczyk}
We disprove this conjecture.
\begin{thm}[store=disprove]\label{thm:disprove}
	There is a class of subcubic trees that has bounded radius-$r$ blind cop-width for every $r\in \mathbb N$ and unbounded pathwidth, thus having bounded blind flip-width and unbounded linear clique-width.
\end{thm} 
We approach this
conjecture by its implications in the weakly sparse setting. Here is a useful proposition due to Gurski~\cite{Gurski2006a}, showing the equivalence of bounded linear clique-width and bounded pathwidth for weakly sparse classes of graphs.
\begin{prop}[Gurski~{\cite[Theorems 9 and 10]{Gurski2006a}}]\label{prop:pw-lcwd}
	Let $G$ be a graph and let $k$ be an integer.
	\begin{enumerate}[label=\rm(\roman*)]
		\item\label{itm:pw-to-lcwd} If $G$ has pathwidth at most $k$, then the linear clique-width of $G$ is at most $k+2$.
		\item\label{itm:lcwd-to-pw} If $G$ has linear clique-width at most~$k$ and has no $K_{n,n}$ subgraph for some $n>1$, then the pathwidth of $G$ is at most $2k(n-1)$.
	\end{enumerate}
\end{prop}

The following lemma shows that bounded blind cop-width implies bounded blind flip-width. This is a blind analog of the inequality between the radius-$r$ flip-width and the radius-$r$ cop-width \cite[Lemma V.5]{T23}.

\begin{lem}\label{bfw-bcw}
	For every $r\in\mathbb{N}\cup\{\infty\}$ and every graph $G$, 
	$\bfw_r(G)\le \bcw_r(G)+2^{\bcw_r(G)}$.
\end{lem}

\begin{proof}
	Given a winning strategy for the cops, one can simulate it using flips in
	the following way; at round $i$, if $k$ cops are at positions $C_i$, given
	the partition made of the singletons of every cop position and
	$V(G)\setminus C_i$ partitioned in at most $2^k$ sets, each consisting of
	vertices having the same set of neighbors in $C_i$, one can perform a flip
	that removes every edge adjacent to a cop without introducing any new edge. 
	In particular, the robber cannot move along a path going through some vertex
	in $C_i\cap C_{i+1}$ before round $i+1$, so he will eventually end up in a
	vertex occupied by a cop, hence isolated, and this strategy is also winning
	for the flipper.
\end{proof}

Thus, \zcref{conj:torunczyk} would imply that any class of graphs with bounded
blind cop-width also has bounded pathwidth, but we give a counterexample to
this property. The counterexample is a sequence of subdivisions of arbitrarily
large complete binary trees (which have unbounded pathwidth by \zcref{pathwidth-tree}) so that three
cops are always sufficient to win the blind cop-width game of radius~$1$.
Here is a theorem of Bernshteyn and Lee~{\cite[Theorem 6]{BL22}} in a weaker form.
\begin{thm}[Bernshteyn and Lee~{\cite[Theorem 6]{BL22}}]\label{thm:tree}
    Every tree has a subdivision whose radius\nobreakdash-$1$ blind cop-width is at most~$3$.
\end{thm}

In particular, we get the following corollary, which will be sufficient for our needs.

\begin{cor}\label{ex:trees_bcw_3}
    For every positive integer $i$, 
    there is a subdivision $T_i$ of a complete binary tree of height $i$ whose radius-$1$ blind cop-width is at most $3$. 
\end{cor}

We give a sketch of its direct proof. A similar idea can be used to prove \zcref{thm:tree}.

\begin{figure}
	\centering
	\begin{tikzpicture}[
		node distance=1cm,
		n/.style={circle, fill=black, inner sep=1.5pt},
		root_node/.style={circle, fill=black, inner sep=1.5pt, label={[label distance=-1mm]above right:$r$}},
		subtree/.style={
			draw,
			fill=gray!20,
			regular polygon,
			regular polygon sides=3,
			minimum size=2.8cm,
			inner sep=0pt
		},
		cop/.style = {circle,fill=blue,inner sep=2.5pt}
	]

		\node[root_node] (r) at (0,0) {};

		\node[subtree] (TL) at (-2,-2) {$T_{i,1}$};
		\node[subtree] (TR) at (4,-4){$T_{i,2}$}; %
		\node [n] at (TL.north){};
		\node [n] at (TR.north){};
		\draw (r) -- (TL.north);
		\draw (r) -- 
		node [pos=0.1,n] {} 
		node [pos=0.2,n] {} 
		node [pos=0.3,n] {} 
		node [pos=0.4,n] {} 
		node [pos=0.5,n] {} 
		node [pos=0.6,n] {} 
		node [pos=0.7,n] {} 
		node [pos=0.8,n] {} 
		node [pos=0.9,n] {} 
		(TR.north);
		\draw[draw=green!50!black] plot [smooth cycle] coordinates {(-3,-2) (-3.2,-3.2) (-.5,-3.2) (-1.8,-1.8) (-2,-1.5)};
		\draw [draw=green!50!black] plot [smooth cycle] coordinates {(1,-1) (3,-3) (2,-5) (6,-5) (5,-3) (2,-.5)};
		\node [cop] at (-2,-1.3) {};
		\node [cop] at (-2.4,-1.5) {};
		\node [cop] at (-1.5,-1.8) {};
	\end{tikzpicture}
	\caption{The state of the game when the cops are cleaning the first
	subtree.}\label{fig}
\end{figure}
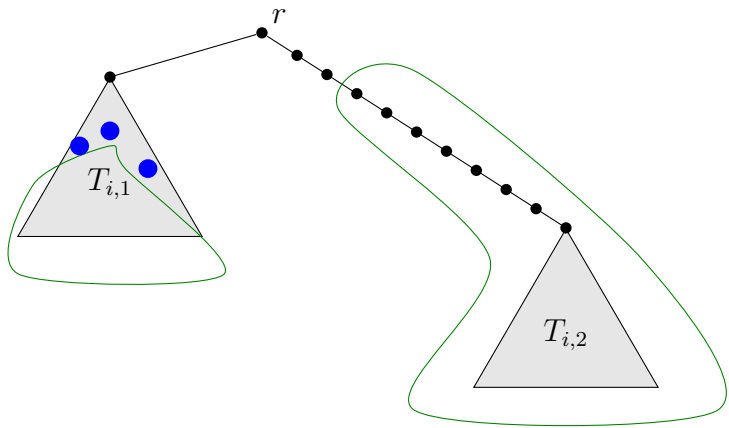

\begin{proof}[Sketch of the Proof]
Let $T_1 = K_1$. For some integer $i\ge 1$, assume that $T_i$ is constructed. Let~$\ell_i$ be the length of a winning strategy for cops in the blind cop-width game of radius 1 on $T_i$.
Given two copies $T_{i,1}, T_{i,2}$ of $T_i$, we build $T_{i+1}$ by creating a new vertex $u$, and
connecting $u$ by an edge to the root of $T_{i,1}$ and by a path $P$ of length $\ell_i$ to the root
of $T_{i,2}$.

The winning strategy for the cops is described as follows.  See \zcref{fig} for an illustration.
\begin{enumerate}
\item One cop stays on $u$, while the two
others clean the path $P$.
\item The three cops recursively clean $T_{i,1}$ top-down. 
\item Two cops clean back $P$,
which was recontaminated while cleaning $T_{i,1}$. 
\item The three cops recursively clean $T_{i,2}$.\qedhere 
\end{enumerate}
\end{proof}

Now let us show how to disprove \zcref{conj:torunczyk}.

\begin{proof}[Proof of \zcref{thm:disprove}]	
	Let $\mathscr{C}$ be the class of subcubic trees having bounded radius-$1$ blind cop-width given by \zcref{ex:trees_bcw_3}.
	By \zcref{feq}, $\mathscr{C}$ has bounded radius-$r$ blind cop-width for every finite $r$.
	By \zcref{pathwidth-tree}, 
	$\mathscr C$ has unbounded pathwidth.

	By \zcref{bfw-bcw}, $\mathscr C$ has bounded blind flip-width
	and 
	by \zcref{prop:pw-lcwd}\ref{itm:lcwd-to-pw}, $\mathscr C$ has unbounded linear clique-width, as trees have no $K_{2,2}$ subgraph.
\end{proof}

Toru\'nczyk \cite[Theorem VI.3 and Corollary IV.6]{T23} proved that the following are equivalent for every class~$\mathscr{C}$ of graphs:
\begin{enumerate}[label=\rm(\roman*)]
\item $\mathscr{C}$ has bounded expansion.
\item There is a function $f:\mathbb{N}\to\mathbb{N}$ such that for every positive integer~$r$
and every graph~$G$ in~$\mathscr{C}$, $\copwidth_r(G)\le f(r)$.
\item $\mathscr{C}$ is weakly sparse and has bounded flip-width. 
\end{enumerate}
We may ask whether the analogous statement regarding the equivalence of (ii) and (iii) is true for blind versions.
In other words, is it true that 
a class of graphs has bounded blind cop-width 
if and only if it is weakly sparse and has bounded blind flip-width?
Actually, this is false as witnessed by the following example.

\begin{prop}\label{bfw_counter_ex}
	For a positive integer $i$, let $T_i$ be a subdivision of a complete binary tree of height~$i$ whose radius-$1$ blind cop-width is at most $3$
	and let $T_i'$ be a tree obtained from $T_i$ by attaching $\abs{V(T_i)}+1$ new vertices to each leaf.
	Then $\{T'_i\}$ has bounded flip-width, while it does not have bounded blind cop-width.
\end{prop}
\begin{proof}%
	Note that $T_i$ exists by \zcref{ex:trees_bcw_3}. 
	By \zcref{pathwidth-tree}, $\pw(T_i)\ge \lceil i/2\rceil$.

	We claim that the radius-$1$ blind flip-width of~$T_i$ is at most $11$.
	Note that a strategy of three cops in the radius-$1$ blind cop-width game can be simulated by $11$-flips, since three vertices induce a partition of the set of other vertices into at most $8$ subsets.
	Once the flipper has restrained the robber to a subtree rooted in a leaf
	of~$T_i$ by simulating the strategy the three cops have for $T_i$, he can
	use one more flip to delete all the new edges adjacent to the leaf of~$T_i$,
	successfully isolating the robber. 
	This means that the radius-$r$ blind flip-width of $\{T_i':i>0\}$ is bounded for every positive integer~$r$.

	On the other hand, we claim that for each $i>0$, 
	$T'_i$ contains a balanced minor model of a complete binary tree of height $\lceil i/2 \rceil$.
	Consider a subtree~$U_i$ of~$T_i$ obtained by, 
	for each internal node~$v$ of odd distance from the root in the initial complete binary tree, 	
	removing the whole subtree rooted at the right child of~$v$
	except for one path from~$v$ to a leaf through the right child. 
	So, $U_i$ contains a minor model $(X_j)$ of a complete binary tree of height~$\lceil i/2\rceil$, where each~$X_j$ contains at least one leaf.
	This minor model can be extended to a balanced minor model in $T_i'$
	because $T_i'$ was obtained by attaching many leaves to each leaf of $T_i$.
	Thus for the function~$g$ in \zcref{tw-tree},
	if $\lceil i/2\rceil \ge g(k)$, then $T_i'$ has radius-$1$ blind cop-width larger than~$k$.
	\zcref{feq} then implies that the radius-$r$ blind cop-width of $\{T'_i:i> 0\}$ is unbounded for each~$r\ge1$.
\end{proof}

\subsection{Blind flip-width of graphs of bounded maximum degree}

We have seen in \zcref{bfw_counter_ex} that bounded blind cop-width is not equivalent to bounded blind flip-width, not only for weakly sparse graphs, but also for $K_{2,t}$-free graphs for $t\ge2$.
However, we can recover functional equivalence if we require the stronger condition that the graphs are $K_{1,t}$-free, or, equivalently, have bounded maximum degree.

\begin{prop}
	If $G$ is a graph of maximum degree $\Delta$, then 
	$\bcw_1(G) \leq 2(\Delta+1)^2 \bfw_3(G)$.
\end{prop}
\begin{proof}
	Let $G$ be a graph of maximum degree at most $\Delta$ and radius-$3$ blind flip-width~$k$.
	Let us show that $2k(\Delta+1)^2$ cops can win on $G$ in the blind cops and robber game.
 	We fix a winning strategy witnessing $\bfw_3(G)\leq k$.
	We may assume that $G$ is connected, because the blind cop-width of a graph is equal to its maximum over all its connected components.

	We employ the following strategy for the cops: at each step, stay on all vertices of flipped parts of size at most $2(\Delta+1)$ and the neighbors of these vertices. There are at most
	$2k(\Delta+1)^2$ such vertices. Let us show that this ensures that the set of possible moves of the robber in the blind cop-width game is a subset
	of those in the radius-$3$ blind flipper game.

	Let $u\rightarrow v$ be a possible move of the robber in the blind cop-width game between rounds $i$ and $i+1$, and
	let $P_u$ and $P_v$ be the parts of $u$ and $v$ in the partition of the $i$-th round in the flipper
	game (note that we may have $P_u=P_v$).
	Since $u\to v$ is a valid move, we deduce that $uv\in E(G)$ and $u$ was not guarded by cops in the $i$-th
	round.
	In particular, given the strategy used in the blind
	cop-width game, both $P_u$ and $P_v$ must contain more than $2(\Delta + 1)$ vertices.
	If the edge	$uv$ is not flipped in the $i$-th round, then $u\to v$ is also a valid move in the $i$-th round in the radius-$3$ flipper game.

	Therefore, we may assume that $(P_u,P_v)$ was flipped in the $i$-th round. 
	Let $G_i$ be the graph obtained after applying the $i$-th round flip to~$G$. First, since 
	\[ 
		\abs{N_{G_i}(u)\cap P_v}= \abs{P_v\setminus (N_G(u)\cup\{u\})} \ge \abs{P_v}-(d_G(u)+1) > 0,
	\]
	there is $w\in P_v$ such that $uw \in E(G_i)$. 
	Similarly, 
	\[ 
		\abs{P_u\cap N_{G_i}(w)\cap N_{G_i}(v)} 
		= \abs{P_u\setminus (N_G(w)\cup N_G(v)\cup\{u,w\})}
		\ge
		|P_u|-(d_G(w)+d_G(v)) > 0,
	\] 
	so there is $t \in P_u$ such that $wt, tv \in E(G_i)$ and $uwtv$ is a
	$u$-$v$ path of length $3$ in $G_i$, which is a valid move in the flipper game of radius $3$.

	Hence, the possible moves of the robber between consecutive rounds in the
	blind cop-width game using this strategy are a subset of the possible moves
	of the robber in the radius-$3$ flipper game using a winning strategy. What
	is left to check is that when the robber gets isolated in the flipper game,
	a cop lands on it in our strategy; it might be isolated but safe, which is
	not sufficient to win in the blind cop-width game.

	We may assume that $G$ has at least two vertices. 
	Let $u$ be an isolated vertex in the flipped graph $G_i$, and $v$ be a neighbor of $u$ in $G$.
	Letting $P_u$ and $P_v$ be the parts of $u$ and~$v$ in the last flip, $(P_u,P_v)$ must then
	be a flipped pair. Now if $|P_v|>2(\Delta+1)$, $u$ would have neighbors in~$P_v$ after flipping,
	which is not possible since it is isolated. Hence $|P_v|\leq 2(\Delta+1)$, 
	so according to our strategy, $u$ is covered by cops, and they win.
\end{proof}

\section{Subdivisions of graphs of bounded treewidth}\label{subdivisions}

\zcref{tw} shows that high blind cop-width can often be explained by
structural reasons, namely high treewidth.
This means that a graph of high
treewidth has high blind cop-width no matter how much we subdivide the
graph, as subdividing preserves treewidth. 
However, by~\zcref{tw-tree}, complete binary trees of high pathwidth also have high blind cop-width, even though they have treewidth $1$.
But this one can be considered less structural, as those trees of high blind
cop-width can always be subdivided to reduce the radius-$1$ blind cop-width to
$3$ by \zcref{ex:trees_bcw_3}. 

In this section, we generalize \zcref{ex:trees_bcw_3}, showing that any graph of treewidth $k$ has a subdivision whose radius-$1$ blind cop-width is at most $k+3$.

Before stating our result, let us recall the concept of a nice tree decomposition
from Bodlaender and Kloks~\cite{BK1996} and Kloks~\cite[Chapter 13]{Kloks1994}.
A tree decomposition $(T,\mathcal X)$ with $\mathcal X=(X_v)_{v\in V(T)}$ of a graph $G$ is a \emph{nice} tree decomposition if $T$ is rooted and the following are true.
\begin{enumerate}[label=\rm (\alph*)]
\item Every node of $T$ has at most $2$ children.
\item If a node $u$ has two children $v$ and $w$, then $X_u = X_v = X_w$. 
\item If a node $u$ has only one child $v$, then either $\abs{X_u}=\abs{X_v}+1$ and $X_v\subseteq X_u$, or $\abs{X_u}=\abs{X_v}-1$ and $X_u\subseteq X_v$.
\end{enumerate}
In a nice tree decomposition~$(T,\mathcal X)$, 
we say a node $u$ of $T$ is a \emph{start} node if it is a leaf,
a \emph{join} node if it has two children,
a \emph{forget} node if it has one child $v$ and $\abs{X_v}>\abs{X_u}$, 
and 
an \emph{introduce} node if it has one child $v$ and $\abs{X_v}<\abs{X_u}$.

\begin{lem}[Kloks~{\cite[Lemma 13.1.2]{Kloks1994}}]\label{lem:nice-tree-decomp}
	Every $n$-vertex graph with tree-width at most $k$ 
	has a nice tree decomposition of width at most $k$ with at most $4n$ nodes.
\end{lem}

To obtain a subdivision of a graph, we need to determine, for each edge~$e$, the length of a path to be replaced with $e$. The following lemma provides such a length, which we will later use to prove the main result of this section.

\begin{lem}\label{lem:subdividing}
	Let $G$ be a graph having a nice tree decomposition $(T,\mathcal X)$ of width at most $k$. Assume that $T$ is an ordered rooted tree, so that we can distinguish the left child and the right child of a join node.
	For a node $u$ of~$T$, let $T[u]$ be the subtree of $T$ induced by all descendants of $u$ and $V_u=\bigcup_{w\in V(T[u])} X_w$.
	For a node $u$ of $T$ and its child $v$, let $E_{u,v}$ be the set of edges  of~$G$ joining a vertex in $X_u$ with a vertex in $V_v\setminus X_u$
	and let $I_{u,v}$ be the set of edges of~$G$ incident with a vertex in $V_v\setminus X_u$. 
	Then there exists a function $\ell:E(G)\to \mathbb N$ such that 
	for every join node $a$ of $T$ with left child $b$ and right child $c$ and for every edge $e\in E_{a,c}$, 
	we have 
	\begin{equation}\label{eq:length}
		\ell(e)\ge 2^{\abs{E_{a,c}}} \abs{V(T[b])}
		\bigl(1+\sum_{f\in I_{a,b}} (\ell(f)-1) \bigr).
	\end{equation}
\end{lem}
\begin{proof}
	Let us define a dependency directed graph $D=(E(G),A)$ on the set of edges of $G$ as follows. 
	For two edges $e$ and $f$, we say $f$ is adjacent to $e$ (equivalently, $(f,e)\in A$) if there is an inequality for $\ell(e)$ in \eqref{eq:length} where the right-hand side has $\ell(f)$.
	Then $(f,e)\in A$ if and only if there exists a join node $a$ with a left child $b$ and a right child $c$ such that $f\in I_{a,b}$ and $e\in E_{a,c}$.

	Now, it is enough to show that $D$ is a directed acyclic graph. 
	Note that $E_{a,b}\subseteq I_{a,b}$.

	First we claim that if an edge $e$ of~$G$ belongs to $I_{a,b}$ for a join node~$a$ of~$T$ with a child~$b$, then every node $u$ of $T$ containing both ends of $e$ is contained in $T[b]$.
	Assume that $e\in I_{a,b}$. Let $x$, $y$ be the ends of $e$ where $x\in V_{b}\setminus X_{a}$. 
	Since $x\in V_{b}\setminus X_{a}$ , every node $v$ of $T$ with $x\in X_v$ belongs to $T[b]$, by the definition of a tree decomposition. 
	This proves the first claim. 

	For an edge $e$, let us write $r(e)$ for the node $u$ closest to the root such that $X_u$ contains both ends of $e$.  
	Let $\preceq$ be a pre-order traversal of $T$. 

	Suppose that $(f,e)\in A$.
	Then there is a join node~$a$ of~$T$ with a left child~$b$ and a right child~$c$ such that $f\in I_{a,b}$ and $e\in E_{a,c}$.
	Then by the claim, 	$r(f)\in V(T[b])$ and $r(e)\in V(T[c])$.
	Since $b$ is the left child and $c$ is the right child of $a$, the traversal $\preceq$ visits $a$, then completely traverses $T[b]$, and then completely traverses $T[c]$. 
	Since $V(T[b])\cap V(T[c])=\emptyset$,
	every node in~$V(T[b])$ comes before every node in $V(T[c])$.
	Therefore we have $r(f) \preceq r(e)$. 

	Suppose that $D$ contains a directed cycle $e_1e_2\cdots e_k e_1$ for $e_1,e_2,\ldots,e_k\in E(G)$. 
	This implies that $r(e_1)\preceq r(e_2)\preceq \cdots \preceq r(e_k)\preceq r(e_1)$. This contradicts the fact that $\preceq$ is a strict total order. 
	This proves that $D$ is acyclic. 

	So there is a topological ordering $\preccurlyeq$ of $V(D)=E(G)$.
	Let $e_1\preccurlyeq e_2\preccurlyeq e_3\preccurlyeq \cdots \preccurlyeq e_m$ be the ordering of $E(G)$ given by $\preccurlyeq$. 
	Let $\mathcal J(e)$ be the set of all triples $(a,b,c)\in V(T)^3$ such that 
	$a$ is a join node, $b$ is a left child of $a$, $c$ is a right child of $a$, and $e\in E_{a,c}$. 
	Thus, for $i=1$ to $m$, we can determine 
	\[ \ell(e_i) = \max \left( \{1\} \cup  \Bigl\{ 2^{\abs{E_{a,c}}}\abs{V(T[b])}  \bigl(1+ \sum_{f\in I_{a,b}} (\ell(f)-1) \bigr) : (a,b,c)\in \mathcal J(e_i)\Bigr\} \right),\]  
	because the value of $\ell(f)$ was already computed by the property that $f\preccurlyeq e_i$ whenever $f\in I_{a,b}$. This completes the proof.
\end{proof}

We are now ready to prove the result. To get some intuition of what is happening
in the proof, one might think of it as applying a similar strategy as the one
of \zcref{thm:disprove} on a subdivision of a graph $G$. The subdivision and the strategy
are based on a nice tree decomposition of $G$.

\begin{thm}[store=treewidthsub]\label{thm:tree-sub}
	Every graph $G$ of treewidth $k$ has a subdivision $H$ whose 
	radius-$1$ blind cop-width is at most $k+3$.
\end{thm}

\begin{proof}
	We may assume that $G$ is connected.
	Let $G$ be a graph of treewidth at most $k$, 
	By \zcref{lem:nice-tree-decomp}, 
	there is a nice tree decomposition $(T,\mathcal X)$ with $\mathcal X=(X_v)_{v\in V(T)}$ of width at most $k$.
	We may further assume that for the root $u$ of $T$, $\abs{X_u}=1$,
	because otherwise we can repeatedly create a parent of the root with one less vertices in its branch set.

	Note that $\abs{X_v}\le k+1$ by the definition of treewidth.
	We may assume that $T$ is an ordered rooted binary tree so that we can distinguish the left child and the right child of a join node.

	Let $\ell:E(G)\to\mathbb N$ be the function given by \zcref{lem:subdividing}.
	Let $H$ be a subdivision of $G$ obtained by replacing each edge~$e$ by a path $P_e$ of length $\ell(e)$.

	For a node $u$ of~$T$, let $T[u]$ be the rooted subtree of $T$ induced by all descendants of~$u$.
	Let $V_u\subseteq V(G)$ be the set $\bigcup_{w\in V(T[u])} X_w$.
	Let $V_u'$ be a subset of $V(H)$ such that 
	\[ V'_u= V_u \cup \bigcup_{e\in E(G[V_u])\setminus E(G[X_u])} V(P_e)	
	.\] 

	We claim that for every $u\in V(T)$, there is a winning strategy $(\mathcal A_u,\mathcal C_u)$
	with $\mathcal A_u=(A_{u,i})_{i\in [m_u]}$ and $\mathcal C_u=(C_{u,i})_{i\in [m_u]}$ 
	for the radius-$1$ blind cop-width game  on $H$ 
	such that 
	\begin{enumerate}[label=\rm (\alph*)]
		\item $A_{u,1}=V'_{u}\setminus C_{u,1}$, $A_{u,i+1}=(A_{u,i}\cup N_H(A_{u,i}))\setminus C_{u,i+1}$ for every $i\in [m_u-1]$ and $A_{u,m_{u}}=\emptyset$, 
		\item $C_{u,1}=X_u$, $\abs{C_{u,i}}\le k+3$, and $A_{u,i}\cap X_u=\emptyset$ for all $i\in [m_u]$, and 
		\item $m_u\le \abs{V(T[u])}(1+\abs{V_u'}-\abs{V_u})$.
	\end{enumerate}
	
	If this claim is true, then by applying it to the root, we deduce that the radius-$1$ blind cop-width of $G$ is at most $k+3$, proving the theorem.

	To prove this claim, we proceed by  induction on $\abs{V(T[u])}$.
	If $u$ is a start node, then $V_u'=X_u$ and therefore $A_{u,1}=\emptyset$ 
	and we are done with $m_u=1$.

	If $u$ is an introduce node, let $v$ be the child of $u$. Then there is a vertex $p\in X_u$ such that $X_u=X_v\cup \{p\}$. 
	By the definition of a tree decomposition, 
	$N_G(p)\cap V_u \subseteq X_v$ and therefore 
	$V_u'=V_v'\cup\{p\}$ 
	and $V_u'\setminus X_u= V_v'\setminus X_v$.
	By the induction hypothesis applied to $v$, there is a winning strategy $(\mathcal A_v,\mathcal C_v)$ satisfying our conditions for $v$
	where $A_{v,1} = V'_v\setminus X_v$. %
	Thus we take  $C_{u,1}:=X_u$
	and $C_{u,i}:=C_{v,i}$ for all $1<i\le m_v$. 
	Then $A_{u,1}=V_u'\setminus X_u = V_v'\setminus X_v$ 
	and by the induction we can see that $A_{u,i+1}=A_{v,i+1}$ for all $i\in [m_v-1]$. 
	Thus, we can take $m_u:=m_v$. 
	Observe that $A_{v,i}\cap X_v=\emptyset$ for all $i\in [m_v]$ and therefore $p\notin A_{u,i}$ for all $i\in [m_u]$. In other words, $p$ is cleaned in
    the beginning, occupied by a cop and cannot be reinfected later.

	If $u$ is a forget node with a unique child $v$,
	then let $(\mathcal A_v, \mathcal C_v)$ be the strategy obtained by the induction hypothesis for $v$
	and let $p$ be the unique vertex in $X_v\setminus X_u$.	
	Note that $V'_u=V'_v \cup \bigcup\{V(P_e):\text{$e$ is an edge of $G[X_v]$ incident with $p$}\}$.
	Let $p_1,p_2,\ldots,p_t$ be a sequence of all vertices in $V'_{u}\setminus V'_v$ so that for every edge~$e$ in $E(G[X_v])$ incident with~$p$, all internal vertices of $P_e$ appear consecutively according to the order in $P_e$ in this sequence. 
	If $t\le 2$, then let $C_{u,1}:=X_u$, $C_{u,2}:=X_v\cup (V'_u\setminus V_v')$,
	and $C_{u,1+i}:=C_{v,i}$ for all $1<i\le m_v$, 
	so that $A_{u,2}=V'_v\setminus X_v$ and $A_{u,m_v+1}=\emptyset$.
	If $t>2$, let $C_{u,1}:=X_u$ and $C_{u,i+1}:=X_v\cup \{p_i,p_{i+1}\}$ for each $i\in [t-1]$.
	Then $A_{u,t}=V_v'\setminus X_v$. 
	We define $A_{u,t-1+i}:=A_{v,i}$ for every $2\le i\le m_v$. 
	Then $A_{u,t-1+m_v}=A_{v,m_v}=\emptyset$.
	In all cases, $m_u\le m_v+ \max(t-1,1) \le \abs{V(T[v])}(1+\abs{V_v'}-\abs{V_v})+ \max(t-1,1) \le \abs{V(T[u])}(1+\abs{V_u'}-\abs{V_u})$.

	Now let us assume that $u$ is a join node with the left child $v$ and the right child $w$.
	Let $E_{u,w}$ be the set of edges of~$G$ joining a vertex in $X_u$ with a vertex in $V_w\setminus X_u$.
	Let $I_{u,v}$ be the set of edges of~$G$ incident with a vertex in $V_v \setminus X_u$.
	Note that \[ \sum_{f\in I_{u,v}} (\ell(f)-1)= \abs{V_v'}-\abs{V_v}.\] 
	For convenience, let $M=\abs{V(T[v])}(1+\sum_{f\in I_{u,v}}(\ell(f)-1))=\abs{V(T[v])}(1+\abs{V_v'}-\abs{V_v})$.
	Recall that by the construction of $H$, every edge $e$ in $E_{u,w}$ has been replaced with a path $P_e$ of length at least $2^{\abs{E_{u,w}}}M$.
	Let $e_1$, $e_2$, $\ldots$, $e_q$ be the edges in~$E_{u,w}$.
	Note that $e_i$ has exactly one end in~$X_u$, say~$x_i$.
	For each $i\in [q]$ and $j\in [2^{q-i}M]$, 
	let $v_{i,j}$ be the $j$-th internal vertex from $x_i$ in~$P_{e_i}$.

	Let $a_0=0$ and $a_i=a_{i-1}+2^{q-i}M$ for $i\in [q]$.
	Since $2^{q-1}+2^{q-2}+\cdots+2^{0} = 2^{q}-1$, we have $a_q=(2^{q}-1)M$.
	Let $v_1,v_2,\ldots,v_{(2^{q}-1)M}$ be the sequence of vertices of $H$ such that 
	for each $i\in [q]$ and $j\in [2^{q-i}M]$, 
	\[ v_{a_{i-1}+j} = v_{i,j}.\]
	Let $(\mathcal A_v,\mathcal C_v)$, $(\mathcal A_w,\mathcal C_w)$ be the winning strategies for the cops obtained by applying the induction hypothesis to $v$ and $w$, respectively. 
	Now, for each $i\in [(2^q-1)M+m_v+m_w]$, let 
	\[
		C_{u,i}:= \begin{cases}
			X_u &\text{if }i=1,\\
			X_u \cup \{v_i,v_{i-1}\} &\text{if }1<i\le (2^q-1)M, \\ 
			X_u\cup \{v_{i-1}\} & \text{if } i = (2^q-1)M+1,\\
			C_{v,i-(2^q-1)M} & \text{if }(2^q-1)M+1<i\le (2^q-1)M+m_v,\\ 
			X_u & \text{if }i = (2^q-1)M+m_v+1,\\
			C_{w,i-(2^q-1)M-m_v} & \text{if }(2^q-1)M+m_v+1<i\le (2^q-1)M+m_v+m_w.
		\end{cases}
	\]

	Now observe that for each $i\in [q]$,
	\[ 
		A_{u,a_i+1} \cap \{v_{i,1},v_{i,2},\ldots,v_{i,2^{q-i}M}\}=\emptyset.
	\]
	As none of $v_{1,1},v_{1,2},\ldots,v_{1,2^{q-i}M}$ are in $C_{u,a_i+j}$
    for any $j>1$, 
	we deduce that for each $1<j\le (2^{q}-1)M-a_i+1$, we have 
	\[ 
		(V'_w\setminus A_{u,a_i+j})\cap \{v_{i,1},v_{i,2},\ldots,v_{i,2^{q-i}M}\}
		= \{ v_{i,1},v_{i,2},\ldots,v_{i,2^{q-i}M-(j-1)}\}.
	\]
	Note that $(2^{q}-1)M-a_i=  \sum_{j=0}^{q-i-1} 2^j M = (2^{q-i}-1)M$.
	Therefore for each $i\in [q]$, 
	\[ 
		(V'_w\setminus A_{u,(2^{q}-1)M+1})\cap\{v_{i,1},v_{i,2},\ldots,v_{i,2^{q-i}M}\} = \{ v_{i,1},v_{i,2},\ldots,v_{i,M}\}.
	\] 
	Thus, $V'_w\setminus A_{u,(2^{q}-1)M+1}=\{v_{i,j}: i\in [q],~j\in [M]\}$.
	In other words, while cleaning the first segments of those $q$ paths, recontamination does occur, but since the speed of the robber is limited, it cannot reach distance $M$ from $X_u$ in $V_v'\setminus X_u$.
    This allows us to buy enough time to clean up $V_v'$ while maintaining $X_u$ clean.
	Thus, while we follow the strategy from $(\mathcal A_v,\mathcal C_v)$, 
	vertices in $X_v$ remain clean, meaning that they are out of $A_{u,i}$ for any $i$ with $(2^q-1)M+1<i\le (2^q-1)M+m_v$.
	Afterwards, by the correctness of $(\mathcal A_v,\mathcal C_v)$, 
	$A_{v,(2^q-1)M}$ does not contain any vertex in $V_v'$.
	Thus, $A_{v,(2^q-1)M+1}\subseteq V_w'\setminus X_w$.
	From that moment we follow the strategy from $(\mathcal A_w,\mathcal C_w)$
	and by the assumption, we end up having $A_{(2^q-1)M + m_v+m_w}=\emptyset$, thus proving that $(\mathcal A_u,\mathcal C_u)$ is a winning game for the cops.

	By the induction hypothesis, we know that $m_v\le \abs{V(T[v])}(1+\abs{V_v'}-\abs{V_v})$ 
	and $m_w\le \abs{V(T[w])}(1+\abs{V_w'}-\abs{V_w})$.
	By definitions and properties of tree decompositions (recall that $X_u=X_v=X_w$), $\abs{V_u'}-\abs{V_u}=\abs{V_v'}-\abs{V_v}+\abs{V_w'}-\abs{V_w}$ and $\abs{V(T[u])}=\abs{V(T[v])}+\abs{V(T[w])}+1$.
	Clearly, $(2^q-1)M \le \abs{V_u'}-\abs{V_u}$, because $v_1,v_2,\ldots,v_{(2^q-1)M}\in V_u'\setminus V_u$.
	Therefore 
	\begin{align*}
		\lefteqn{
		\abs{V(T[u])}(1+\abs{V_u'}-\abs{V_u})
		}\\
		&= (\abs{V(T[v])}+\abs{V(T[w])}) (1+\abs{V_v'}-\abs{V_v}+\abs{V_w'}-\abs{V_w}) + 1+\abs{V'_u}-\abs{V_u}\\
		&\ge m_v+m_w + (2^q-1)M.
	\end{align*}
	Therefore $m_u:=(2^q-1)M+m_v+m_w\le \abs{V(T[u])}(1+\abs{V_v'}-\abs{V_v})$, as required.
	This proves the claim. This completes the proof.
\end{proof}

\section{Connections to other games and parameters}\label{connections}

\paragraph{Inspection number.}
Bernshteyn and Lee~\cite{BL22} introduced 
the \emph{zero-visibility search game} and the inspection number. 
In this game, cops and robber alternate turns. 
On the cops' turn, they search a set of at most $k$ vertices, and
win if the robber is present in this set.
Otherwise, on the robber's turn,
the robber can move along any edge. The minimum number $k$ such that cops
have a winning strategy on a graph $G$ is the \emph{inspection number} of~$G$, denoted by~$\IN(G)$. 
The inspection number is called the \emph{firefighter number} by Althoetmar, Schade, and Schürenberg~\cite{ASS2025}.

Now, let us recall a precise formulation from \cite{BL22}.
A \emph{$k$-search} of length~$\ell$ on a graph~$G$ is a sequence $\mathcal{S}=(S_1,S_2,\ldots,S_\ell)$ of subsets of~$V(G)$ where $\abs{S_i}\le k$ for all $i\in[\ell]$.
Given a $k$-search $\mathcal{S}=(S_1,S_2,\ldots,S_\ell)$ of subsets of $V(G)$
of size at most $k$, a sequence $(\fc_i)_{i=0}^\ell$ of sets of \emph{fully cleared vertices} and a sequence $(\pc_i)_{i=0}^\ell$ of sets of \emph{pre-cleared vertices} are defined recursively as follows: 
\begin{align*}
	\fc_0(\mathcal S) &= \emptyset,\\
	\pc_{i}(\mathcal S) &= \fc_{i-1} \cup ~S_{i} &\text{ for }i\in [\ell],\\
	\fc_{i}(\mathcal S) &= \pc_{i}(\mathcal S)\setminus~N_G(V(G)\setminus \pc_{i}(\mathcal S)) &\text{ for }i\in [\ell].
\end{align*}
A $k$-search~$\mathcal S$ of length~$\ell$ is called \emph{successful} if $\fc_\ell(\mathcal S) = V(G)$, which implies $\pc_\ell(\mathcal S) = \pc_{\ell+1}(\mathcal S)= V(G)$.

\begin{prop}[Bernshteyn and Lee~{\cite[Proposition 10]{BL22}}]\label{prop:inspection}
The inspection number of a graph~$G$ is equal to the minimum number~$k$ such
that there exists a successful $k$-search in $G$.
\end{prop}

This game might seem different, but the fact that
cops abandon their position after their turn actually makes them equivalent
to the radius-$1$ blind cop-width.
\begin{prop}\label{prop:bcw_in}
	For every graph~$G$, $\IN(G) = \bcw_1(G)$.
\end{prop}

\begin{proof}
    First, let us reformulate the radius-$1$ blind cop-width. 
	If $(\mathcal A,\mathcal C)$ is a strategy for the cops in the radius-$1$ blind cop-width, then it is easy to see that 
    \eqref{eq:bcw_r} is equivalent to the following with the same initial condition $A_0=V(G)$:
    \[A_{i} = \left(A_{i-1} \cup N_G(A_{i-1})\right) \setminus C_{i} \text{ for all $i\ge 1$.}\]
	The cops win if $A_\ell=\emptyset$ for some $\ell$.

	Suppose that $\mathcal C=(C_i)_{i=1}^\ell$ is a sequence of subsets of $V(G)$. 
	Now we claim that 
	$\pc_{i}(\mathcal C)= V(G)\setminus A_{i}$ for all $i\ge 1$.
	We proceed by induction on~$i$.
	Clearly, it is true if $i=1$. 
	Now let us assume that $i>1$. Then, by the induction hypothesis, 
	\[ \fc_{i-1}(\mathcal C)= (V(G)\setminus A_{i-1}) \setminus N_G(A_{i-1}) 
	= V(G)\setminus (A_{i-1}\cup N_G(A_{i-1})).
	\] 
	Thus, we deduce that 
	\[ 
		\pc_{i}(\mathcal C) = V(G)\setminus (A_{i-1}\cup N_G(A_{i-1}))\cup C_i 
		= V(G)\setminus ((A_{i-1}\cup N_G(A_{i-1}) \setminus C_{i} ))
		= V(G)\setminus A_{i}.
	\] 
	
	If $A_\ell=\emptyset$, then $\pc_{\ell}(\mathcal C)=V(G)$ and therefore $\fc_{\ell}(\mathcal C)=V(G)$.
	Conversely, if $\fc_{\ell}(\mathcal C)=V(G)$, then 
	$\pc_{\ell}(\mathcal C)=V(G)$ and therefore $A_{\ell}=\emptyset$. 
	Thus $\IN(G)=\bcw_1(G)$.
\end{proof}

We define the \emph{topological blind cop-width} of a graph~$G$, denoted by $\tbcw(G)$, also known as the \emph{topological inspection number} in \cite{BL22}.
For a graph~$G$, $\tbcw(G)$ is the minimum
$k$ such that for every subdivision $H$ of~$G$, there is a further subdivision
$H'$ of~$H$ with $\bcw_1(H')\leq k$.

Bernshteyn and Lee~{\cite[Proposition 23]{BL22}} showed that 
$\tbcw(G)\le \abs{V(G)}+2$. This is now a corollary of \zcref{thm:tree-sub}, as $\tbcw(G)\le \tw(G)+3\le (\abs{V(G)}-1)+3$.

Here are two of the problems proposed by Bernshteyn and Lee~\cite{BL22}.
\bla*
\blb* 
\zcref{tw} answers \zcref{prob68} positively and \zcref{prob69} negatively. 
Indeed, as treewidth is invariant under subdivisions, \zcref{thm:tree-sub,tw}
prove the following.
\begin{prop}
    The topological blind cop-width is functionally equivalent to treewidth. 
\end{prop}

For complete graphs, Bernshteyn and Lee found the following upper bound.

\begin{prop}[Bernshteyn and Lee~{\cite[Proposition 67]{BL22}}]
\label{prop:kn}
    For $n\ge 1$, 
    $\tbcw(K_n)\le \lceil n/3\rceil + 2$.
\end{prop}

Bernshteyn and Lee~{\cite[Theorem 32]{BL22}} showed that $\tbcw(K_4)\ge 4$ 
by a rather technical proof. 
The following proposition immediately implies that $\tbcw(K_n)=\Theta(n)$, and also answers \zcref{prob68}.
\begin{prop}\label{prop:hadwiger}
    $\lceil \frac{n+6}{4}\rceil \le \tbcw(K_n) \le \lceil \frac{n+6}{3}\rceil$ for all integers $n\ge 2$.
\end{prop}

\begin{proof}%
    \zcref{prop:kn} produces a subdivision $H'$ of~$K_n$ with $\bcw_1(H') \le \lceil n/3\rceil + 2$.
    This gives the upper bound.
	For the lower bound, we may assume that $n\ge 7$ because $\tbcw(K_2)\ge 2$ and $\tbcw(K_3)\ge 3$.
    \zcref{lem:balancedcliqueminor} together with
    \zcref{prop:hbtobcw} gives the lower bound.
\end{proof}

Since \zcref{prop:hadwiger}  shows that the
number of cops needed to clean $K_n$ if the cops are allowed to choose a
subdivision to begin with is between $n/4$ and $n/3+3$, one may wonder if the cop
number of $\tw(G)+3$ given by \zcref{thm:tree-sub} could be optimized.

\zcref{prop:ktree} shows that there are graphs having treewidth $k$ and
topological blind cop-width at least $(k+3)/2$,
showing that the complete graph on $k+1$ vertices is not the graph of treewidth $k$
with the maximum topological blind cop-width.

We remark that despite the fact that both \zcref{lem:balancedcliqueminor,prop:hbtobcw} are close to optimal, we suspect that the upper bound of \zcref{prop:kn} is closer to the truth.

\paragraph{Zero-visibility cop number.}
A common variant of cops and robber games, particularly in
the game theory community, considers alternating turns and restricts
cops to move along edges with speed~$1$ as well. 

The associated parameter denoted by $c_0$ is called the \emph{zero-visibility cop number} 
by Dereniowski, Dyer, Tifenbach, and Yang~\cite{DDTY15}, and is also studied in~\cite{XYZZ19,XYZ22}.
It is not functionally equivalent to blind cop-width, as it is implied by~\cite[Theorem 3.8]{DDTY15b}
that trees of high pathwidth have high zero-visibility cop number, 
while \zcref{ex:trees_bcw_3} 
provides trees of radius-$1$ blind cop-width at most $3$ and arbitrarily large pathwidth.
This example is in fact present in~\cite{DDTY15}, as the addition of a
universal vertex can be seen as giving cops infinite speed at the cost of
doubling that of the robber, which does not change much by \zcref{feq}. 

\begin{prop}\label{prop:c0}
	For every graph~$G$, $\bcw_1(G) \leq 2 c_0(G)$.
\end{prop}
The factor of $2$ is necessary, as unlike in the zero-visibility search game, here cops stay at vertices when the robber is trying to move.
Notably, $\IN(K_2)=\bcw_1(K_2)=2$ while $c_0(K_2)=1$.

\begin{proof}
Notice that in the blind setting, whether cops announce their moves or
not makes no difference since we always consider all the possible strategies of the
robber, so we might forget about this condition. The restriction on
the movement only makes it harder for cops in the zero-visibility
cops and robber game, so we can also forget about it to prove this
inequality.

Given a winning strategy $(S_i)_{i=1}^\ell$ in the zero-visibility 
cops and robber game, we can set $C_i := S_i\cup S_{i+1}$ in the
radius-1 blind cop-width game.
This way, in the helicopter game the robber will be prevented from 
moving through cops originally at step $i$ in the zero-visibility game 
as they would not be lifted, so at most $2c_0(G)$ cops will win.
\end{proof}

As a corollary of \zcref{tw} and \zcref{prop:c0}, we get the following lower bound on $c_0$.

\begin{cor}[store=lbcnot]
For every integer $k$, there is an integer $h(k)$ such that every graph of treewidth
	at least $h(k)$ has zero-visibility cop number larger than~$k$.  
\end{cor}

In order to obtain an upper bound on $c_0$ in terms of blind cop-width, we need to restrict ourselves
to graphs of bounded diameter.

\begin{prop}
For a graph~$G$ of diameter at most $d$, we have $c_0(G)\le 4d \bcw_1(G)$.
\end{prop}

\begin{proof}
	First let us show that $c_0(G)\le 2\bcw_d(G)$.
	Let $k=\bcw_d(G)$.
Consider a winning strategy $(S_i)_{i=1}^\ell$ for $k$ cops in the radius-$d$ blind cop-width game, and let $(A_i)_{i=1}^\ell$
be the sequence of possible positions of the robber in a game using this strategy. We may assume that $S_i\neq\emptyset$ for all $i$.
Here is a strategy for the cops in the zero-visibility cop number game. 
Start by putting~$k$ cops on $S_1$ and $k$ cops on $S_2$. 
Then for each integer $i\in[\ell]\setminus\{1,2\}$ and rounds $d(i-3)+2$ to $d(i-2)+1$, make $k$ cops stay on $S_{i-1}$, and another $k$ cops move from their current positions to $S_{i}$
following some paths of length at most $d$ in $G$. Letting $(B_i)_{i=1}^{d(\ell-2)+1}$ be the sequence of
possible positions of the robber in this game, we get that $B_{d(i-2)+1}\subseteq A_i$, in particular $B_{d(\ell-2)+1}=\emptyset$.
Thus, $c_0(G)\le 2\bcw_d(G)$.
Using the explicit bound $\bcw_d(G) \leq 2d \bcw_1(G)$ from \zcref{feq}, we get the desired result.
\end{proof}

\paragraph{Hunting number.}
We end this section by mentioning another interesting variation of the zero-visibility search game, studied by Abramovskaya, Fomin, Golovach, and Pilipczuk~\cite{AFGP16} under the name \emph{Hunters and Rabbit game}.
The game is defined as the zero-visibility search game, with the additional condition
that the robber is forced to move at each turn. 
Since this condition can only make it easier for cops to win, the associated \emph{hunting number}, denoted by $\HH$, clearly satisfies 
$\HH(G)\leq \IN(G)=\bcw_1(G)$. 
On the other hand, the trees of \zcref{bfw_counter_ex} are an example of graphs with high blind
cop-width, hence inspection number, for which four cops suffice in the Hunters
and Rabbit game, so we do not get functional equivalence.

\begin{prop}
    For a positive integer $i$, let $T_i$ be a subdivision of a complete binary tree of height~$i$ whose radius-$1$ blind cop-width is at most $3$
	and let $T_i'$ be a tree obtained from $T_i$ by attaching $\abs{V(T_i)}+1$ new vertices to each leaf.
    Then $\{T_i'\}$ has bounded hunting number, while it does not have bounded blind cop-width.
\end{prop}

\begin{proof}
    	Note that $T_i$ exists by \zcref{ex:trees_bcw_3}.
    The fact that $\{T_i'\}$ has unbounded blind cop-width follows from \zcref{bfw_counter_ex}.

    We claim that the hunting number of $T_i$ is at most 4. A winning strategy is
    given by letting 3 cops simulate a winning strategy for $T_i$. Each time 
    the robber is restrained to $T'_i[u]\setminus \{u\}$ for some leaf $u$
    (when $u$
    would be cleaned in the blind cop-width game), the fourth cop then stays on $u$
    in the next round (recall that $T[u]$ denotes the subtree of $T$ rooted at
    node $u$). That way, if the robber was hiding in a child of $u$, it would have
    to move and be shot by the fourth cop in the next round.
\end{proof}

We do however retrieve functional equivalence
when restricted to bounded degree graphs.

\begin{thm}\label{maxdeghunt}
    For any graph $G$ without isolated vertices, $\bcw_1(G)\leq (\Delta(G)+1)\HH(G)$.
\end{thm}
\begin{proof}
Let $C_1,\dots, C_t$ be a winning strategy for the hunters and rabbit game.
We will show that $N[C_1],\dots,N[C_t]$ is a winning cop strategy for the blind cop-width game of radius $1$.
To do this, we inductively show that the set of possible rabbit locations in the hunters and rabbits game at step~$i$ is a superset of the set of possible robber locations given this new strategy.
The base case is trivial.

Suppose for contradiction that at step $i+1$ there is some possible robber location $v$ that is not a possible rabbit location. 
By induction, every possible robber location at step $i$ is a possible rabbit location at step $i$, so $v$ is only a possible location for a robber at step $i+1$ by staying still.
That is to say, $v$ is a possible robber/rabbit location at step~$i$ that is not adjacent to another possible robber/rabbit location at step~$i$.
Let $a_1,\dots a_i$ be a possible path for the robber ending at $v$, so that $a_i=v$ and each $a_j$ for $j<i$ is a possible robber location at step~$j$ which is equal to or adjacent to $a_{j+1}$.
Let $i'$ be the smallest integer such that $a_{i'}$ is not adjacent to another possible rabbit location at step $i'$.
If $i'=1$, then $C_1$ contains all neighbors of $a_{1}$, so $N[C_1]$ contains $a_{1}$, contradicting that $a_1$ is a possible robber location at step $1$.
Otherwise, by induction both $a_{i'-1}$ and some neighbor $w$ of $a_{i'-1}$ are possible rabbit locations at step $i'-1$.
Since $a_{i'-1}$ is not a possible rabbit location at step $i'$ (by definition of $i'$, as it is a neighbor of $a_{i'}$), we must then have $a_{i'-1}\in C_{i'}$, as otherwise a rabbit could go from $w$ to $a_{i'-1}$. Hence $a_{i'}\in N[C_{i'}]$, which again is a contradiction.

Thus the set of possible robber locations at each step is a subset of the possible rabbit locations, meaning that since $C_1,\dots, C_t$ is a winning strategy for the hunters $N[C_1],\dots , N[C_t]$ is a winning strategy for the cops.
Note that $\abs{N[C_i]}\leq (\Delta+1)\abs{C_i}$, so this completes the proof.
\end{proof}

We have the following corollary of \zcref{tw}, generalizing the results of Abramovskaya, Fomin, Golovach, and Pilipczuk~\cite[Theorem 2]{AFGP16} on
the hunting number of grids.
\begin{cor}[store=lbhunt]
    For every integer $k$, there is an integer~$j(k)$ such that every graph of treewidth
	at least~$j(k)$ has hunting number larger than~$k$. 
\end{cor}

\begin{proof}
    Let $g$ be the function of \zcref{tw}. Let $k$ be a positive integer.
	By \zcref{grid}, 
	there is an integer $j$ such that every graph~$G$ of treewidth at least~$j$ contains a subdivision~$H$ of a wall as a subgraph and $\tw(H)\ge g(4k)$.
	By \zcref{tw}, $\bcw_1(H) > 4k$, but $\Delta(H) = 3$, so by \zcref{maxdeghunt}, $\HH(H)> k$. Since the hunting number is monotone under taking subgraphs~\cite{AFGP16}, we have $\HH(G)>k$.
\end{proof}

 Dissaux, Fioravantes, Galhawat, and Nisse~\cite[Theorem~6]{DFGN25} proved that every tree has a subdivision whose hunting number is at most $2$. As a corollary of \zcref{thm:tree-sub}, we get the following.

\begin{cor}[store=huntsub]
	Every graph $G$ of treewidth $k$ has a subdivision $H$ whose 
	hunting number is at most $k+3$.
\end{cor}

\section{Discussions}\label{sec:end}

The results of \zcref{connections} and \zcref{tw} allow us to obtain
the comparisons of parameters.
Indeed, by \zcref{feq}, for every $r\in \mathbb{N}$, $\bcw_r$ is functionally equivalent to $\bcw_1$ and so we may omit $\bcw_r$ for any finite $r$. 
By \zcref{tw}, bounded radius-$1$ blind cop-width implies bounded treewidth.
Also note that \zcref{ex:trees_bcw_3} separates $\bcw_1$ from $\pw$.
Thus we have the following. 
\[ 
	\copwidth_1  \leq  \copwidth_2 \leq \cdots < \tw \leq \bcw_1 < \pw.
\] 
Here, the relation $\leq$ between graph parameters should be understood as $p\leq q$ if there exists some function $f:\mathbb{N}\to \mathbb{N}$ such that $p(G)\leq f(q(G))$ for all graphs $G$.

\zcref{tw} shows that the blind cop-width of a graph lies asymptotically
somewhere between its treewidth and pathwidth, and can be pushed to either extreme by
subdividing edges appropriately. 
In terms of algorithmic complexity, the class of problems which are tractable on bounded pathwidth graphs is strictly broader than the class of problems which are tractable on bounded treewidth graphs (see \cite{BKLMO22, B25}).
The same may be true of graphs with bounded blind cop-width, although we suspect not in any meaningful sense given how subdivisions affect the parameter.
We are more interested in understanding the structural properties of the parameter, which prompt several intriguing questions for future work.

In the direction of upper bounds, it would be interesting to give a better description of the
winning strategies for the cops. Currently, the only general upper bound we have is in terms of the
pathwidth, which is obtained by considering a winning strategy in the case the robber could
move at infinite speed. The strategy is then described by a path decomposition of the graph, as
the cops win by successively occupying the vertices in the branch sets of the decomposition.
Generally speaking, winning strategies are often associated with decompositions of the graphs.
If few cops can win, then the graph has some structure and can be decomposed in a way that
follows a winning strategy.
A naive attempt to describe an appropriate decomposition of the blind cop-width game would be to say that a sequence of vertex sets is a $\bcw$-decomposition if and only if it corresponds to a winning strategy. Note that the analogous definition for the pathwidth game does not correspond to the notion of path decompositions. Path decompositions, in general, correspond to winning strategies with some
specific minimality conditions.
As a starting point then, it would be helpful to develop an understanding of what ``minimal'' winning strategies for the blind cop-width game look like, which we are currently lacking.

\begin{problem}
    What is the structure of graphs of radius-$1$ blind cop-width at most~$k$?
\end{problem}

In the dual direction, it would be very interesting to better understand lower bounds as well.
Given that we have shown some obstructions to low blind cop-width, it is natural to ask if they
are exactly the right ones.

\begin{problem}
    Is there an unbounded function $f$ such that each graph of blind cop-width at least~$k$
    contains a binary tree of height at least $f(k)$ as a balanced minor?
\end{problem}

It seems that simple cases may provide good intuitions
on this problem, so we recall the following question that we only answered in
an asymptotic way.

\pkn*
As pointed out in \zcref{prop:hadwiger}, we have a lower bound of $\lceil (n+6)/4\rceil$ and an upper bound of $\lceil (n+6)/3\rceil $ for $\tbcw(K_n)$.

Finally, while we have shown that topological blind cop-width is functionally equivalent to
treewidth, it could be interesting to get some more precise characterizations, at least up to
constant terms.

\bibliography{bib}

\providecommand{\bysame}{\leavevmode\hbox to3em{\hrulefill}\thinspace}
\providecommand{\MR}{\relax\ifhmode\unskip\space\fi MR }
\providecommand{\MRhref}[2]{%
  \href{http://www.ams.org/mathscinet-getitem?mr=#1}{#2}
}
\providecommand{\href}[2]{#2}
\begin{thebibliography}{10}

\bibitem{AFGP16}
Tatjana~V. Abramovskaya, Fedor~V. Fomin, Petr~A. Golovach, and Micha{\l} Pilipczuk, \emph{How to hunt an invisible rabbit on a graph}, European J. Combin. \textbf{52} (2016), no.~part A, 12--26. \MR{3425961}

\bibitem{ASS2025}
Julius Althoetmar, Jamico Schade, and Torben Sch{\"u}renberg, \emph{Complexity of firefighting on graphs}, arXiv:2505.11082, 05 2025.

\bibitem{BKLMO22}
R\'{e}my Belmonte, Eun~Jung Kim, Michael Lampis, Valia Mitsou, and Yota Otachi, \emph{Grundy {D}istinguishes {T}reewidth from {P}athwidth}, SIAM J. Discrete Math. \textbf{36} (2022), no.~3, 1761--1787. \MR{4457712}

\bibitem{BL22}
Anton Bernshteyn and Eugene Lee, \emph{Searching for an intruder on graphs and their subdivisions}, Electron. J. Combin. \textbf{29} (2022), no.~3, Paper No. 3.9, 46. \MR{4446719}

\bibitem{B25}
Hans~L. Bodlaender, Carla Groenland, Hugo Jacob, Lars Jaffke, and Paloma~T. Lima, \emph{{XNLP}-completeness for parameterized problems on graphs with a linear structure}, Algorithmica \textbf{87} (2025), no.~4, 465--506.

\bibitem{BK1996}
Hans~L. Bodlaender and Ton Kloks, \emph{Efficient and constructive algorithms for the pathwidth and treewidth of graphs}, J. Algorithms \textbf{21} (1996), no.~2, 358--402. \MR{98g:68122}

\bibitem{BG19}
Jessalyn Bolkema and Corbin Groothuis, \emph{Hunting rabbits on the hypercube}, Discrete Mathematics \textbf{342} (2019), no.~2, 360--372.

\bibitem{CT2020}
Julia Chuzhoy and Zihan Tan, \emph{Towards tight(er) bounds for the excluded grid theorem}, J. Combin. Theory Ser. B \textbf{146} (2021), 219--265. \MR{4155282}

\bibitem{DDTY15b}
Dariusz Dereniowski, Danny Dyer, Ryan~M. Tifenbach, and Boting Yang, \emph{The complexity of zero-visibility cops and robber}, Theoret. Comput. Sci. \textbf{607} (2015), no.~part 2, 135--148. \MR{3426926}

\bibitem{DDTY15}
\bysame, \emph{Zero-visibility cops and robber and the pathwidth of a graph}, J. Comb. Optim. \textbf{29} (2015), no.~3, 541--564. \MR{3316705}

\bibitem{DiestelGT}
Reinhard Diestel, \emph{Graph theory, 4th edition}, Graduate texts in mathematics, vol. 173, Springer, 2012.

\bibitem{DFGN23}
Thomas Dissaux, Foivos Fioravantes, Harmender Gahlawat, and Nicolas Nisse, \emph{Recontamination helps a lot to hunt a rabbit}, MFCS 2023--48th International Symposium on Mathematical Foundations of Computer Science, Schloss Dagstuhl-Leibniz-Zentrum f{\"u}r Informatik, 2023, pp.~42--1.

\bibitem{DFGN25}
Thomas Dissaux, Foivos Fioravantes, Harmender Galhawat, and Nicolas Nisse, \emph{Further results on the hunters and rabbit game through monotonicity}, Information and Computation (2025), 105302.

\bibitem{EST1994}
John~A. Ellis, Ivan~H. Sudborough, and Jonathan~S. Turner, \emph{The vertex separation and search number of a graph}, Inform. and Comput. \textbf{113} (1994), no.~1, 50--79. \MR{1283019 (95d:05115)}

\bibitem{GH2023}
Cyril Gavoille and Claire Hilaire, \emph{Minor-universal graph for graphs on surfaces}, arXiv:2305.06673, 05 2023.

\bibitem{Gurski2006a}
Frank Gurski, \emph{Linear layouts measuring neighbourhoods in graphs}, Discrete Math. \textbf{306} (2006), no.~15, 1637--1650. \MR{2251098}

\bibitem{HMV2004}
Darrel Hankerson, Alfred Menezes, and Scott Vanstone, \emph{Guide to elliptic curve cryptography}, Springer Professional Computing, Springer-Verlag, New York, 2004. \MR{2054891}

\bibitem{Kinnersley1992}
Nancy~G. Kinnersley, \emph{The vertex separation number of a graph equals its path-width}, Inform. Process. Lett. \textbf{42} (1992), no.~6, 345--350. \MR{1178214 (93d:05143)}

\bibitem{KP85}
Lefteris~M. Kirousis and Christos~H. Papadimitriou, \emph{Interval graphs and searching}, Discrete Math. \textbf{55} (1985), no.~2, 181--184. \MR{798534}

\bibitem{Kloks1994}
Ton Kloks, \emph{Treewidth}, Lecture Notes in Computer Science, vol. 842, Springer-Verlag, Berlin, 1994, Computations and approximations. \MR{1312164 (96d:05038)}

\bibitem{NO12}
Jaroslav Ne{\v{s}}et{\v{r}}il and Patrice Ossona~de Mendez, \emph{Sparsity - graphs, structures, and algorithms}, Algorithms and combinatorics, vol.~28, Springer, 2012.

\bibitem{Reitwiesner1960}
George~W. Reitwiesner, \emph{Binary arithmetic}, Advances in {C}omputers, {V}ol. 1, Academic Press, New York-London, 1960, pp.~231--308. \MR{122018}

\bibitem{grid}
Neil Robertson and Paul Seymour, \emph{Graph minors. {V}. {E}xcluding a planar graph}, J. Combin. Theory Ser. B \textbf{41} (1986), no.~1, 92--114. \MR{89m:05070}

\bibitem{RST1994}
Neil Robertson, Paul Seymour, and Robin Thomas, \emph{Quickly excluding a planar graph}, J. Combin. Theory Ser. B \textbf{62} (1994), no.~2, 323--348. \MR{96c:05050}

\bibitem{Scheffler1989}
Petra Scheffler, \emph{Die {B}aumweite von {G}raphen als ein {M}a\ss\ f\"ur die {K}ompliziertheit algorithmischer {P}robleme}, Report MATH, vol. 89-04, Akademie der Wissenschaften der DDR, Karl-Weierstrass-Institut f\"ur Mathematik, Berlin, 1989, With an English summary. \MR{1004243}

\bibitem{ST93}
Paul Seymour and Robin Thomas, \emph{Graph searching and a min-max theorem for tree-width}, J. Combin. Theory Ser. B \textbf{58} (1993), no.~1, 22--33. \MR{94b:05197}

\bibitem{TU15}
Satoshi Tayu and Shuichi Ueno, \emph{On evasion games on graphs}, Discrete and Computational Geometry and Graphs - 18th Japan Conference, {JCDCGG} 2015, vol. 9943, 2015, pp.~253--264.

\bibitem{T23}
Szymon Toru{\'n}czyk, \emph{Flip-width: Cops and robber on dense graphs}, 2023 IEEE 64th Annual Symposium on Foundations of Computer Science (FOCS), 2023, pp.~663--700.

\bibitem{XYZZ19}
Yuan Xue, Boting Yang, Farong Zhong, and Sandra Zilles, \emph{A partition approach to lower bounds for zero-visibility cops and robber}, Combinatorial algorithms, Lecture Notes in Comput. Sci., vol. 11638, Springer, Cham, 2019, pp.~442--454. \MR{3991255}

\bibitem{XYZ22}
Yuan Xue, Boting Yang, and Sandra Zilles, \emph{A simple method for proving lower bounds in the zero-visibility cops and robber game}, J. Comb. Optim. \textbf{43} (2022), no.~5, 1545--1570. \MR{4452885}

\bibitem{XYZ24}
Yuan Xue, Boting Yang, and Sandra Zilles, \emph{The zero-visibility cops and robber game on graph products}, Theoretical Computer Science \textbf{1007} (2024), 114676.

\end{thebibliography}
\bibliographystyle{amsplain}

\end{document}